\documentclass{amsart}
\usepackage{amsmath,amssymb,amscd,verbatim}

\newcommand{\Z}{\ensuremath{\mathbb{Z}}}
\newcommand{\U}{\ensuremath{U(\mathfrak{g})}}
\newcommand{\Uh}{\ensuremath{U_{h}(\mathfrak{g}) } }
\newcommand{\Uqt}{\ensuremath{U_{h}^{qt}(\mathfrak{g}) } }
\newcommand{\Uhp}{\ensuremath{U_{h}(\mathfrak{g_{+}}) } }
\newcommand{\Uhb}{\ensuremath{U_{h}(\mathfrak{b_{+}}) } }
\newcommand{\Uhbhat}{\ensuremath{U_{h}(\tilde{\mathfrak{b}}_{+}) } }

\newcommand{\Hom}{\ensuremath{\operatorname{Hom}}}
\newcommand{\End}{\ensuremath{\operatorname{End}}}
\newcommand{\g}{\ensuremath{\mathfrak{g}}}
\newcommand{\ghat}{\ensuremath{\tilde{\mathfrak{g}}}}
\newcommand{\liea}{\ensuremath{\mathfrak{a}}}
\newcommand{\h}{\ensuremath{\mathfrak{h}}}
\newcommand{\nilp}{\ensuremath{\mathfrak{n}}}
\newcommand{\borel}{\ensuremath{\mathfrak{b}}}
\newcommand{\borelp}{\ensuremath{\mathfrak{b}_{+}}}
\newcommand{\borelhat}{\ensuremath{\tilde{ \mathfrak{b}}}}
\newcommand{\gl}{\ensuremath{\mathfrak{gl}(m|n)}}
\newcommand{\Pf}{\ensuremath{{}^{\prime}\text{f}}}
\newcommand{\f}{\text{f}}
\newcommand{\pu}{\ensuremath{{}^{\prime}\text{U}}}
\newcommand{\cu}{\text{U}}

\newcommand{\C}{\ensuremath{\mathbb{C}} }

\newcommand{\N}{\ensuremath{\mathbb{N}} }
\newcommand{\p}[1]{\ensuremath{\bar {#1}}}
\newcommand{\M}{\ensuremath{\mathcal{M_{\g}}}}
\newcommand{\A}{\ensuremath{\mathcal{A}}}
\newcommand{\D}[1]{\ensuremath{\Delta({#1})} }

\newcommand{\Uhatp}{\mathcal{\widetilde{U}_{+}}}
\newcommand{\Up}{\mathcal{U_{+}}}
\newcommand{\Upqd}{\mathcal{U_{+}^{\star}}}
\newcommand{\cas}{\Omega }
\newcommand{\rmap}{\operatorname{p} }
\newcommand{\rmapim}{\text{Im}\rmap}
\newcommand{\comp}[1]{\widehat{#1}}
\newcommand{	\qd}{{}^{\star}}
\newcommand{	\qdo}{{}^{\star op}}

\newcommand{\DJg}{\ensuremath{U_{h}^{DJ}(\g) } }
\newcommand{\DJb}{\ensuremath{U_{h}^{DJ}(\borel_{+}) } }

\newcommand{\ttp}{\hat{\otimes}}
\newcommand{\atp}{\tilde{\otimes}}
\newcommand{\kzt}{A_{\g,t}}
\newcommand{\Mt}{\ensuremath{\mathcal{M}_{\g}^{t}}}
\newcommand{\At}{\ensuremath{\mathcal{A}^{t}}}
\newcommand{\flip}{\operatorname{\tau}}
\newcommand{\s}{\ensuremath{\flip}}
\newcommand{\tF}{\widetilde{F}}
\newcommand{\tMt}{\ensuremath{\widetilde{\mathcal{M}}_{\g}^{t}}}

\newcommand{\kz}{A_{\g}}

\newcommand{\OH}{{H}}
\newcommand{\OUhp}{{U}_{h}(\g_{+})}
\newcommand{\OUh}{{U}_{h}}

\newcommand{\OF}{{F}}
\newcommand{\OFq}{\overline{F}}

\newcommand{\OJ}{{J}}
\newcommand{\lform}{{C}}
\newcommand{\Fiso}{{\Psi}}
\newcommand{\Otheta}{\vartheta}
\newcommand{\FH}{\overline{H}}
\newcommand{\FUhp}{\overline{U}_{h}(\g_{+})}
\newcommand{\FUh}{\overline{U}_{h}}
\newcommand{\FF}{\overline{F}}
\newcommand{\FJ}{\overline{J}}
\newcommand{\Fvtheta}{\overline{\vartheta}}
\newcommand{\Ftheta}{\vartheta}
 \newtheorem{thm}{Theorem}
 \newtheorem{cor}[thm]{Corollary}
\newtheorem{lem}[thm]{Lemma}

\newtheorem{prop}[thm]{Proposition}

\theoremstyle{definition}
\newtheorem{defn}[thm]{Definition}

\newtheorem{rem}[thm]{Remark}

\begin{document}

 \title{Etingof-Kazhdan quantization of Lie superbialgebras}
 \author{Nathan Geer} \address{School of Mathematics,  
Georgia Institute of Technology, 
Atlanta, GA 30332-0160} \email{geer@math.gatech.edu}

\begin{abstract}
For every semi-simple Lie algebra $\g$ one can construct the Drinfeld-Jimbo algebra $\DJg$.  This algebra is a deformation Hopf algebra defined by generators and relations.  To study the representation theory of $\DJg$, Drinfeld used the KZ-equations to construct a quasi-Hopf algebra $\kz$.  He proved that particular categories of modules over the algebras $\DJg$ and $\kz$ are tensor equivalent.  Analogous constructions of the algebras $\DJg$ and $\kz$ exist in the case when $\g$ is a Lie superalgebra of type A-G.  However, Drinfeld's proof of the above equivalence of categories does not generalize to Lie superalgebras.  In this paper, we will discuss an alternate proof for Lie superalgebras of type A-G.  Our proof utilizes the Etingof-Kazhdan quantization of Lie (super)bialgebras.  It should be mentioned that the above equivalence is very useful.  For example, it has been used in knot theory to relate quantum group invariants and the Kontsevich integral.
\end{abstract}

\maketitle

\section{Introduction}  %some words here:  WRITE ABOUT THE IDEAS, NOT THE OBJECTS.
Quantum groups were introduced independently by Drinfeld and Jimbo around 1984.  One of the most important examples of quantum groups are deformations of universal enveloping algebras.  These deformations are closely related to Lie bialgebras.  In particular, every deformation of a universal enveloping algebra induces a Lie bialgebra structure on the underling Lie algebra.  In \cite{D7} Drinfeld asked if the converse of this statement holds:  ``Does there exist a universal quantization for Lie bialgebras?''  Etingof and Kazhdan gave a positive answer to this question.  In this paper we further this work by extending Etingof and Kazhdan's work from Lie bialgebras to the setting of Lie superbialgebras.  Moreover, we will generalize a theorem of Drinfeld's from Lie algebras to Lie superalgebras of type A-G.  

\vspace{10pt}

Given a semisimple Lie algebra $\g$, Drinfeld \cite{D5} constructs the following algebras:
\begin{enumerate}
\item the Drinfeld-Jimbo quantization $\DJg$ of $\g$ which is a deformation of the Hopf algebra $U(\g)$, 
\item a quasi-Hopf algebra $\kz$ which is isomorphic as a vector space to $U(\g)[[h]]$.   
\end{enumerate}
These algebras are quite different in nature.  $\DJg$ is a Hopf algebra which is defined algebraically by generators and relations.  The non-trivial and complicated structure of $\DJg$ is encoded in these relations and the formulas defining the coproduct.   On the other hand, the definition of $\kz$ is based on the theory of the Knizhnik-Zamolodchikov differential equations.  Compared to $\DJg$, the algebra structure and the coproduct of $\kz$ are easy to define.   The rich structure of $\kz$ is encoded in the fact that its coproduct is not coassociative. 

Drinfeld was interested in the representation theory of the algebras $\DJg$ and $\kz$.  Let $X$ be a topological algebra and let $X \text{-} Mod_{fr}$ be the category of topologically free $X$-modules of finite rank (see \ref{SS:TopoFreeMod}).  

\begin{thm}[\cite{D5}]\label{T:Drinfeld}
The categories $\DJg  \text{-} Mod_{fr}$ and $\kz  \text{-} Mod_{fr}$ are tensor equivalent.
\end{thm}

This theorem allows one to play the differences of $\DJg$ and $\kz$ off of one another, leading to a deeper understanding of the category $\DJg  \text{-} Mod_{fr}$.  It turns out that  Theorem \ref{T:Drinfeld} is also useful in knot theory.  In particular, Le and Murakami used Theorem \ref{T:Drinfeld} to show that quantum group knot invariants arising from representations of Lie algebras can be studied through the Kontsevich integral. 

The algebras $\DJg$ and $\kz$ can be constructed for some classes of Lie superalgebras.  In \S \ref{S:KZ} we will construct $\kz$ for a suitable Lie superalgebra $\g$.  The generalization of $\DJg$ to the setting of Lie superalgebras has been considered by many authors (see \cite{FLV,KTol,Yam91}).  This generalization introduces defining relations (e.g. (\ref{E:QserreC}-\ref{E:QserreD})) that are of a different form than the standard quantum Serre relations of $\DJg$.   Let us call these additional relations the extra quantum Serre-type relations.  Unlike the case for semi-simple Lie algebras, the (quantum) Serre-type relations are not well understood for all Lie superalgebras.  For this reason we will consider the Lie superalgebras of type A-G.  Yamen \cite{Yam91,Yam94} obtained (quantum) Serre-type for every Lie superalgebra of type A-G.   

The proof of Theorem \ref{T:Drinfeld} does not have a straightforward generalization to the setting of Lie superalgebras. 
Drinfeld's proof uses deformation theoretic arguments based on the fact that $H^{i}(\g,U(\g))=0, \; i=1,2,$ for a semisimple Lie algebra.  In general, this vanishing result is not true for Lie superalgebras (for example $\mathfrak{sl}(2|1)$, \cite{SZ}).  However, in \S \ref{S:ModuleCat} we will prove that Thoerem \ref{T:Drinfeld} is true when $\g$ is a Lie superalgebra of type A-G.   Our proof is based on a different approach than Drinfeld's, utilizing the quantization of Lie (super)bialgebras.    

Our proof of Thoerem \ref{T:Drinfeld} (when $\g$ is a Lie superalgebra of type A-G) starts by generalizing the Etingof-Kazhdan quantization of Lie bialgebras to the setting of Lie superbialgebras.  Note that it can be shown that $\g$ can be given a natural structure of a Lie superbialgebra.  Let $U_{h}(\g)$ be the E-K quantization of $\g$.  By construction $\Uh$ is gauge equivalent to $\kz$.  With sufficient hypotheses, Drinfeld showed that if two algebras are gauge equivalent then their module categories are equivalent.  In a similar fashion, we will show that $\Uh  \text{-} Mod_{fr}$ is tensor equivalent to $\kz  \text{-} Mod_{fr}$.  As one would expect the generalizations discussed in this paragraph are straightforward.

The proof is completed by constructing a Hopf algebra isomorphism between $U_{h}(\g)$ and $\DJg$.  This method is similar to \cite{EK6} where it is shown that analogous result holds for any generalized Kac-Moody Lie algebra $\liea$.  The proof of \cite{EK6} shows that $U_{h}(\liea)$ is given by generators and relations.  In particular, the authors of \cite{EK6} define a bilinear form and use results of Lusztig \cite{L} to show that the quantum Serre-type relations are in the kernel of this form.  Similar techniques apply in the case when $\g$ is a Lie superalgebra of type A-G.  However, as mentioned above, $\DJg$ has extra quantum Serre-type relations.  In order to adapt the above methods we will extend results of Lusztig \cite{L} to the setting of superalgebras and check directly that the extra quantum Serre-type relations are in the kernel of the appropriate bilinear form.   

The Etingof-Kazhdan quantization \cite{EK1} has two important properties that we use in this paper: the first being that it is functorial and second that it commutes with taking the double.  
With this in mind, we will next discuss the notion of the double of an object.  Let $\g$ be a finite dimensional Lie superbialgebras.   The double of $\g$ is the direct sum $D(\g):=\g\oplus \g^{*}$ with a natural structure of a (quasitriangular) Lie superbialgebras.  Similarly, let $A$ be a quantized universal enveloping (QUE) superalgebra and let $A\qd$ be its quantum dual, i.e. a QUE superalgebra which is dual (in an appropriate sense) to $A$.  The double of $A$ of is the tensor product $D(A):=A \otimes A\qd$ with a natural structure of a quasitriangular QUE superalgebra.  By saying the E-K quantization commutes with taking the double we mean that $D(U_{h}(\g)) \cong U_{h}(D(\g))$ as quasitriangular QUE superalgebras.

\vspace{10pt}

We will now give an outline of this paper.  There are several different quantization given in this paper which turn out to be isomorphic.  We hope that following outline will help the reader understand why each quantization is important.   

In \S \ref{Prelim}, we will recall facts and definitions related to Lie superbialgebras, topologically free modules and QUE superalgebra.  In \S \ref{S:QUg},  we will give the definition of a Lie superalgebra of type A-G and its associated D-J type quantization $\DJg$.  In \S \ref{S:KZ}, we will use the super KZ equations to define the quasi-Hopf superalgebra $\kz$.  We will also define the Drinfeld category.  

In \S\ref{S:NQof}, we will extend the Etingof-Kazhdan quantization of finite dimensional Lie bialgebras, given in Part I of \cite{EK1}, to the setting of Lie superbialgebras.  Let $\g$ be a finite dimensional Lie superbialgebra.  Section \ref{S:NQof} consists of three important parts: (1) the construction of a quantization $\OH$ of the double $D(\g)$, (2) show that $\OH$ has a Hopf sub-superalgebra $\OUh(\g)$ which is a quantization of $\g$,  (3) prove that $\OH$ and $\OUh(\g)$ are further related by the following isomorphism of quasitriangular Hopf superalgebras  
\begin{equation}
\label{E:IntroUhTensorDual}
\OH \cong \OUh(\g) \otimes \OUh(\g)\qd
\end{equation}
where $\OUh(\g)\qd$ is the quantum dual of $\OUh(\g)$ and $D(\OUh(\g)):=\OUh(\g) \otimes \OUh(\g)\qd$.  

In \S \ref{S:QUofquasi}, we will construct the E-K quantization of quasitriangular Lie superbialgebras.  Let us denote this quantization by $\Uqt$ where $\g$ is a quasitriangular Lie superbialgebras.  The quantization $\Uqt$ is similar to the quantization of Lie superbialgebra of section \ref{S:NQof}.   In particular, by construction 
\begin{equation}
\label{E:InrtoUhqtDouble}
U_{h}^{qt}(D(\g)) = \OH
\end{equation}
 for any finite dimensional Lie superbialgebra $\g$. 

In \S \ref{S:2QUofg}, we will construct a second quantization of finite dimensional Lie superbialgebras.  The importance of this quantization is that it is functorial.  It turns out that it is isomorphic to the quantization given in section \ref{S:NQof}.  For this reason we also denote it by $\OUh(\g)$.  

In \S \ref{S:FunctQ}, we will use the functoriality of the quantization to show that $\Uqt \cong \OUh(\g)$ for any finite dimensional quasitriangular Lie superbialgebra $\g$.  As noted above, the double of a Lie superbialgebra has a natural structure of a quasitriangular Lie superbialgebra.  Therefore, we have
\begin{equation}
\label{E:IntroUqtUh}
U_{h}^{qt}(D(\g))  \cong \OUh(D(\g))
\end{equation} 
for any finite dimensional Lie superbialgebras $\g$.  We will close section \ref{S:FunctQ} by combining (\ref{E:IntroUhTensorDual}), (\ref{E:InrtoUhqtDouble}) and (\ref{E:IntroUqtUh}) to conclude that the E-K quantization commutes with taking the double.

In \S \ref{S:QofglConclution}, we will prove that the E-K quantization $U_{h}(\g)$ is isomorphic to the D-J type quantization $\DJg$, where $\g$ is a Lie superalgebra of type A-G.  The proof of this will rely heavily on the fact that the E-K quantization is functorial and commutes with taking the double. 

In \S \ref{S:ModuleCat}, we give a proof of Theorem \ref{T:Drinfeld} when $\g$ is a Lie superalgebra of type A-G.

\subsection*{Acknowledgments} The author is grateful to P. Etingof for his useful suggestions and conversations.  Also he thanks A. Berenstein, B. Shelton and A. Vaintrob for their attention and comments.

\section{Preliminaries}\label{Prelim}
Let $k$ be a field of characteristic zero.  
 
\subsection{Superspaces and Lie super(bi)algebras}\label{S:Lbialg}
In this subsection we recall facts and definitions related to superspaces and Lie super(bi)algebras, for more details see \cite{K,Manin}.

A \emph{superspace} is a $\Z_{2}$-graded vector space $V=V_{\p 0}\oplus V_{\p 1}$ over $k$.  We denote the parity of a homogeneous element $x\in V$ by $\p x\in \Z_{2}$.  We say $x$ is even (odd) if $x\in V_{\p 0}$ (resp. $x\in V_{\p 1}$).  Let $V$ and $W$ be superspaces.  The space of linear morphisms $\Hom_{k}(V,W)$ from $V$ to $W$ has a natural $\Z_{2}$-grading given by $f\in \Hom_{k}(V,W)_{\p j}$ if $f(V_{i})\subseteq W_{i+j}$ for $\p i, \p j \in \Z_{2}$.  In particular, the dual space $V^{*}=\Hom_{k}(V,k)$ is a vector superspace where $k$ is the one-dimensional superspace concentrated in degree $\p0$, i.e. $k=k_{\p 0}$.  Throughout this paper the tensor product will have the natural induced $Z_{2}$-grading.  Let $\flip_{V,W}:V\otimes W \rightarrow W \otimes V $ be the linear map given by
\begin{equation}
\label{E:Flip}
 \flip_{V,W}(v\otimes w)=(-1)^{\p v \p w}w\otimes v
\end{equation}
for homogeneous $v \in V$ and $w \in W$.  When it is clear what $V$ and $W$ are will write $\flip$ for $\flip_{V,W}$.
   A linear morphism can be defined on homogeneous elements and then extended by linearity.  When it is clear and appropriate we will assume elements are homogeneous.    Throughout, all modules will be $\Z_{2}$-graded modules, i.e. module structures which preserve the $\Z_{2}$-grading (see \cite{K}).

%Roughly speaking, the super analogue of an object, whose underlying structure is a vector space, is superspace such that one picks up a negative sign in any formula when two odd elements pass by one another.   
A \emph{Lie superalgebra} is a superspace $\g=\g_{\p 0} \oplus \g_{\p 1}$ with a superbracket $[\: , ] :\g^{\otimes 2} \rightarrow \g$ that preserves the  $\Z_{2}$-grading, is super-antisymmetric ($[x,y]=-(-1)^{\p x \p y}[y,x]$), and satisfies the super-Jacobi identity (see \cite{K}).  
A \emph{Lie superbialgebra} is a Lie superalgebra $\g$ with a linear map $\delta : \g \rightarrow \wedge^{2}\g$ that preserves the $Z_{2}$-grading and satisfies both the super-coJacobi identity and cocycle condition (see \cite{A}).
A triple $(\g,\g_{+},\g_{-})$ of finite dimensional Lie superalgebras is a finite dimensional \emph{super Manin triple} if $\g$ has a non-degenerate super-symmetric invariant bilinear form $<,>$, such that $\g \cong \g_{+}\oplus \g_{-}$ as superspaces, and $\g_{+}$ and $\g_{-}$ are isotropic Lie sub-superalgebras of $\g$.  % i.e. $<x,y>=0$ for all $x,y\in \g_{\pm}$
There is a one-to-one correspondence between finite dimensional super Manin triple and finite dimensional Lie superbialgebra (see \cite[Proposition 1]{A}).

% Quasitriagular Lie superbialgebras
Let $\g$ be a Lie superalgebra. %Let $\s: \U\otimes \U \rightarrow \U\otimes \U$ be the map given by 
%$$x\otimes y  \mapsto (-1)^{\bar{x}  \bar{y} }y\otimes x$$
%for homogeneous $x,y \in \U$.  
Let $r\in \g \otimes \g$ and let
$$CYB(r):=[r_{12},r_{13}]+[r_{12},r_{23}]+[r_{13},r_{23}]\in \g^{3}$$
be the classical Yang-Baxter element.  
A \emph{quasitriangular Lie superbialgebra} is a triple $(\g, [\: ,] , r)$ where $(\g,[\: ,])$ is a Lie superalgebra and $r$ is an even element of $\g \otimes \g$ such that $r+\flip(r)$ is $\g$-invariant, $CYB(r)=0$ and $(\g, [\:,], \partial r)$ is a Lie superbialgebra, where $\partial r(x):=[x\otimes 1 + 1 \otimes x,r]$.

Now we define the double of a finite dimensional Lie superbialgebra.  Let $(\g_{+}, [\: , ]_{\g_{+}}, \delta)$ be a finite dimensional Lie superbialgebra and $(\g,\g_{+},\g_{-})$ its corresponding super Manin triple.  Then $\g:=\g_{+}\oplus \g_{-}$ has a natural structure of a quasitriangular Lie superbialgebra as follows.   The bracket on $\g$ is given by
\begin{align}
\label{R:BracketDouble}
   [x,y] &=	\begin{cases}
		[x,y]_{\g_{+}} & \text{if } x, y \in \g_{+}\\
		[x,y]_{\g_{-}} & \text{if } x, y \in \g_{-}\\
		(ad^{*}x)y - (-1)^{\p x \p y}(1\otimes y)\delta(x) &\text{if } x\in \g_{+},y\in \g_{-}
	\end{cases}
\end{align}
where $ad^{*}$ is the coadjoint action of $\g_{+}$ on $\g_{-}\cong\g_{+}^{*}$.  Let $p_{1},...,p_{n}$ be a homogeneous basis of $\g_{+}$.  Let $m_{1},...,m_{n}$ be the basis of $\g_{-}$ which is dual to  $p_{1},...,p_{n}$, i.e. $<m_{i},p_{j}>=\delta_{i,j}$.  Define $r = \sum p_{i}\otimes  m_{i} \in \g_{+}\otimes \g_{-}\subset \g \otimes \g$.  Then the triple $(\g, [,],  r)$ is a quasitriangular Lie superbialgebra (see \cite{A}).  We call $\g$ the \emph{double} of $\g_{+}$ and denote it by $D(\g_{+})$.

We can also define the Casimir element of $\g$.  Notice $m_{1},...,m_{n},p_{1},...,p_{n}$ is a basis of $\g$ that is dual to the basis $p_{1},...,p_{n},(-1)^{\p m_{1}}m_{1},...,(-1)^{\p m_{n}}m_{n}$.  Define the \emph{Casimir element} to be 
\begin{equation}
\label{D:Cas}
\cas =\sum p_{i}\otimes m_{i}+ \sum (-1)^{\p m_{i}}m_{i}\otimes p_{i}=r + \flip(r).
\end{equation}  An element $a \in \g \otimes \g$ is invariant (resp. super-symmetric) if $[x\otimes 1 + 1 \otimes x, a]=0$ for all $x\in \g$ (resp. $a=\flip(a)$).  The element $\cas$ is an even, invariant, super-symmetric element.  Also, note that the element $\cas$ is independent of the choice of basis $p_{1},...,p_{n}$.

%     subsection        Topologically free modules

\subsection{Topologically free modules} \label{SS:TopoFreeMod} Here we recall the notion of topologically free modules (for more detail see \cite{Kas,ES}).

Let $K=k[[h]]$, where $h$ is an indeterminate and we view $K$ as a superspace concentrated in degree $\p 0$.  %A $K$-module is a superspace $M$ which is a $K$-module such that the action of $K$ preserves the grading.
  Let M be a module over $K$.  Consider the inverse system of $K$-modules
$$p_{n}:M_{n}=M/h^{n}M \rightarrow M_{n-1}=M/h^{n-1}M.$$
Let $\comp{M}=\underleftarrow{\lim} M_{n}$ be the inverse limit.  Then $\comp{M}$ has the natural inverse limit topology (called the $h$-adic topology).  We call $\comp M $ the $h$-adic completion of M.  

%  Def of Topol. Free

Let $V$ be a $k$-superspace.  Let $V[[h]]$ to be the set of formal power series.  The superspace $V[[h]]$ is naturally a $K$-module and has a norm given by
$$|| v_{n}h^{n } + v_{n+1}h^{n+1}+ \cdots ||=2^{-n} \: \: \text{ where } v_{n} \neq 0.$$
The topology defined by this norm is complete and coincides with the h-adic topology.  We say that a $K$-module $M$ is \emph{topologically free} if it is isomorphic to $V[[h]]$ for some $k$-module $V$.  Notice that if $f:M\rightarrow N$ is a $K$-linear map between topologically free modules then $f$ is continuous in the $h$-adic topology since $f(h^{n}M)\subseteq h^{n}N$ by $K$ linearity.  For this reason we will assume that all $K$-linear maps are continuous.  

%  Def of Topol tenor product

Let $M,N$ be topologically free $K$-modules.  We define the \emph{topological tensor product} of $M$ and $N$ to be $\comp{M \otimes_{K} N}$ which we denote by $M\otimes N$.  This definition gives us the convenient fact that $M \otimes N$ is topologically free and that
$$V[[h]]\otimes W[[h]] = (V\otimes W)[[h]]$$
for $k$-module $V$ and $W$.

%Let $\A$ be the symmetric tensor category of topologically free $k[[h]]$-modules, with super commutativity isomorphism and trivial associativity isomorphism.

We say a (Hopf) superalgebra defined over $K$ is topologically free if it is topologically free as a $K$-module and the tensor product is the above topological tensor product.

% Subsection       Super Quantized Universal Enveloping Algebras
 
\subsection{Quantized Universal Enveloping  Superalgebras} A \emph{quantized universal enveloping} (QUE) superalgebra $A$ is a topologically free Hopf superalgebra over $\C[[h]]$ such that $A/hA$ is isomorphic as a Hopf superalgebra to $U(\g)$ for some Lie superalgebra $\g$.  The follow proposition was first given in the non-super case by Drinfeld \cite{D} and latter proven in the super case by Andruskiewitsch \cite{A}.

\begin{prop}[\cite{D},\cite{A}]\label{P:QUE_Lie_bi}
Let $A$ be a QUE superalgebra: $A/hA\cong \U$.  Then the Lie superalgebra $\g$ has a natural structure of a Lie superbialgebra defined by
\begin{equation}
\label{R:ConditionDelta}
 \delta (x)  = h^{-1}(\Delta({\tilde{x}}) - \Delta^{op}(\tilde{x})) \mod h,  \hspace{15pt} x \in \g
\end{equation}
where $\tilde{x}\in A$ is a preimage of $x$ and $\Delta^{op}:=\flip_{\U,\U} \circ \Delta$ (for the definition of $\flip$, see (\ref{E:Flip})).
\end{prop} 

\begin{defn}
Let $A$ be a QUE superalgebra and let $(\g,[,], \delta)$ be the Lie superbialgebra defined in Proposition \ref{P:QUE_Lie_bi}.  We say that $A$ is a quantization of the Lie superbialgebra $\g$.
\end{defn}
%   quasitriangular quantization of $(\g,r)$

Let $A$ be a Hopf superalgebra and let $R\in A \otimes A$ be an invertible homogeneous element.  We say $(A,R)$ is a \emph{quasitriangular} Hopf superalgebra if 
\begin{align}
\label{E:Rmatrix}
   R \Delta =& \Delta^{op} R,   
\end{align}
\begin{align}
\label{R:comultR}
    (\Delta \otimes 1)(R) =& R_{13}R_{23},  & (1 \otimes \Delta )(R)=& R_{13}R_{12}.
\end{align}
From relations (\ref{R:comultR}) it follows that $(\epsilon \otimes 1)R = (1 \otimes \epsilon )R=1$ which implies that $R$ is even.

Let $A$ be a quantization of a quasitriangular Lie superbialgebra $(\g,r)$ and let $R\in A \otimes A$ be an invertible homogeneous element.  We say $(A,R)$ is a \emph{quasitriangular quantization} of $(\g,r)$ if $R$ satisfies (\ref{E:Rmatrix}), (\ref{R:comultR}) and 
\begin{align}
\label{E:RmatModh2}
  R \equiv & 1 + hr \mod h^{2}.    
\end{align}
%note from general principles it follows that $R$ is even.  See intro to thesis.

\subsection{The quantum dual and the double}\label{SS:DualDouble}  In this subsection we define the quantum dual and double of a QUE superalgebra.  We will use these construction throughout the rest of the text.  The definition of the quantum dual was first given by Drinfeld \cite{D} in the non-super case.  For more on quantum duals see \cite{D,ES,Gav}. 

Let $\A$ be the symmetric tensor category of topologically free $k[[h]]$-modules, with the super commutativity isomorphism $\flip$ given in (\ref{E:Flip})
 and the canonical associativity isomorphism.
Let $A$ be a QUE superalgebra and set $A^{*}=\Hom_{\A}(A,k[[h]])$.  Then $A^{*}$ is a topological Hopf superalgebra where the multiplication, unit, coproduct, counit, and antipode are given by $f g(x)=(f \otimes g) \Delta (x)$, $\epsilon$, $\Delta f (x \otimes y)=f(xy)$, $1$, and $S^{*}$ (respectively) for $f,g \in A^{*}$ and $x,y \in A$.  Let $I^{*}$ be the maximal ideal of $A^{*}$ defined by the kernel of the linear map $A^{*}\rightarrow k$ given by $f \mapsto f(1) \mod h$.  This gives a topology on $A^{*}$ where $\{ (I^{*})^{n}, n\geq 0 \}$ is a basis of the neighborhoods of zero.  

%  Def   of  quantum dual

Here we give the definition of the quantum dual.  
Define $(A^{*})^{\vee}$ to be the $h$-adic completion of the $k[[h]]$-module $\sum_{n \geq 0} h^{-n}(I^{*})^{n}$.  Then $(A^{*})^{\vee}$ is a QUE superalgebra which denote by $A\qd$.  We call $A\qd$ the \emph{quantum dual} of $A$.   
Let $\delta_{n} : A \rightarrow A^{\otimes n}$ be the linear map given by $\delta_{1}(a)= a - \epsilon(a) 1$, $\delta_{2}(a)=\Delta (a) -a \otimes 1 - 1 \otimes a + \epsilon (a) 1 \otimes 1$, etc.  Define $A^{\prime}=\{ a \in A | \delta_{n}(a) \in h^{n}A^{n}\}$.  Then as shown in \cite{Gav} we have 
\begin{align}
\label{E:QDualandDual}
   (A\qd)^{\prime}=& A^{*}, &  (A^{\prime})^{\vee}= & A. 
 \end{align}

Now we define the notion of the double of $A$.  
Let $\{x_{i}\}_{i \in I}$ be a basis of $A$ and let $ \{y_{i}\}_{i \in I}$ be the corresponding dual elements of $A^{*}$, i.e. $<y_{i},x_{j}>=\delta_{ij}$.  From \cite[\S 3.5]{Gav} it follows that $\hat R=\sum_{i \in I}x_{i} \otimes y_{i}$ is a well defined element of $A \otimes A\qd$.   

% Def of Quant Double.
The following proposition was first due to Drinfeld.
\begin{prop}\label{P:DoubleQUE}
Let $A$ be a QUE superalgebra and $A\qdo$ its dual QUE superalgebra with opposite coproduct ($\Delta^{op}=\flip_{A,A} \circ \Delta$).  Let $\hat R$ be the canonical element defined above.  Then there exist a unique Hopf superalgebra structure on $D(A):=A \otimes  A\qdo$ such that 
\begin{enumerate}
  \item $A$ and $A\qdo$ Hopf sub-superalgebra of $D(A)$.
  \item The linear map $A \otimes  A\qdo \rightarrow D(A)$ given by $a \otimes a^{\prime} \mapsto aa^{\prime}$ is a bijection. 
  \item $\hat R$ is a quasitriangular structure for $D(A)$. 
\end{enumerate}
\end{prop}
\begin{pf}
The proof follows as in the pictorial proof of Proposition 12.1 in \cite{ES}.  One only needs to notice that the corresponding pictures hold in the super case and account for the necessary signs in relation (12.4) and in the proof of Lemma 12.1.
\end{pf}
We call $D(A)$ the \emph{quantum double} of $A$.

%  marker
\section{The Drinfeld-Jimbo type quantization of Lie superalgebras of type A-G}\label{S:QUg}
In this section we recall the defining relations of both a complex Lie superalgebra of type A-G and it quantum analogue.  The relations defining these superalgebras are not easily obtained and are of a different nature than relations arising from Lie algebras.  In this section we also show that a Lie superalgebras of type A-G has a natural structure of a Lie superbialgebra.  For the purposes of this paper Lie superalgebras of type A-G will be complex and include the Lie superalgebra $D(2,1,\alpha)$.

Any two Borel subalgebras of a semisimple Lie algebra are conjugate.  It follows that semisimple Lie algebras are determined by their root systems or equivalently their Dynkin diagrams.   However, not all Borel subalgebras of classical Lie superalgebras are conjugate.  As shown by Kac \cite{K} a Lie superalgebra can have more than one Dynkin diagram depending on the choice of Borel.  However, using Dynkin diagrams Kac gave a characterzation of Lie superalgebras of type A-G.  Using the standard Borel sub-superalgebra, Floreanini, Leites and Vinet \cite{FLV} were able to construct defining relations for some Lie superalgebras and their quantum analogues.  Then  
Yamane \cite{Yam94} gave defining relations for each Dynkin diagram of a Lie superalgebra of type A-G.  These relations are given by formulas which depend directly on the choice of Dynkin diagram.  For this reason, we will restrict our attention to the simplest case and only consider root systems with at most one odd root.

\subsection{Lie superalgebras of type A-G}\label{SS:g}   
Let $\g:=\g_{\p 0} \oplus \g_{\p 1}$ be a Lie superalgebra of type A-G such that $\g_{\p 1}\neq \emptyset$.  As mentioned above, Kac \cite{K} showed that $\g$ is characterized by its associated Dynkin diagrams or equivalently Cartan matrices.  A Cartan matrix associated to a Lie superalgebra is a pair consisting of a matrix $M$ and a set $\tau$ determining the parity of the generators.  As shown by Kac \cite{K}, there exist simple root systems of $\g$ with exactly one odd root.  Let $\Phi=\{\alpha_{1},...,\alpha_{s}\}$ be such a simple root system and let $(A,\{m\})$ be its corresponding Cartan matrix where $\alpha_{m}$ is the unique odd root.  Note that all simple root systems with exactly one odd root are equivalent and lead to the same Cartan matrix (see \cite[\S 2.5.4]{K}).   The Dynkin diagrams corresponding to such Cartan matrices are listed in Table VI of \cite{K}.  For notational convenience we set $I=\{1,...,s\}$.

\begin{thm}[\cite{Yam91,Yam94}]\label{T:GenRelforLieSuper}
Let $\g$ be a Lie superalgebra of type A-G with associated Cartan matrix $(A=(a_{ij}),\tau)$ where $\tau=\{m\}$ (as above) or $\tau=\emptyset$ (purely even case).  Then $\g$ is generated by $h_{i}$, $e_{i},$ and $f_{i}$ for $i\in I$ (whose parities are all even except for $e_{t}$ and $f_{t}$, $t\in \tau$, which are odd) where the generators satisfy the relations
\begin{align}
\label{R:LieSuperalg}
   [h_{i}, h_{j}]&=0,   & [h_{i}, e_{j}] &=a_{ij}e_{j},   &
   [ h_{i},f_{j}] &=-a_{ij} f_{j} &   [e_{i},f_{j}] &=\delta_{ij}h_{i}
\end{align}
and the ``\emph{super classical Serre-type}'' relations
$$[e_{i},e_{i}] =[f_{i},f_{i}] =0 \text{ for } i\in \tau$$
$$(\text{ad } e_{i})^{1+|a_{ij}|}e_{j}=(\text{ad } f_{i})^{1+|a_{ij}|}f_{j}=0, \text{ if } i \neq j,\text{ and } i \notin \tau $$
\begin{multline*}
[e_{m},[e_{m-1},[e_{m},e_{m+1}]]]=[f_{m},[f_{m-1},[f_{m},f_{m+1}]]]=0 \\ \text{ if $ m-1,m, m+1\in I $ and $  a_{mm}=0$,}
\end{multline*}
%if $ m-1,m, m+1\in I $ and $  a_{mm}=0$,
\begin{multline}
\label{E:ClassicSerre}
[[[e_{m-1},e_{m}] e_{m}],e_{m}]=[[[f_{m-1},f_{m}] f_{m}],f_{m}]=0  \\ \text{ if the Cartan Matrix $A$ is of type B, $\tau = \{m\}$ and $s=m$}.
\end{multline}
where $[,]$ is the super bracket, i.e.
$[x,y]=xy-(-1)^{\p x \p y}yx$.  
\end{thm}

For the rest of this paper when considering Lie superalgebras of type A-G we will assume that these Lie superalgebras are defined by the generators and relations given in Theorem \ref{T:GenRelforLieSuper}.

%      Subsection     $\gl$ as a Lie superbialgebra

\subsection{Lie superbialgebra structure}\label{SS:glbialg} 
In this subsection we will show that Lie superalgebras of type A-G have a natural Lie superbialgebra structure.  The following results are straight forward generalizations of the non-super case. 

Let $\g$ be a  Lie superalgebra of type A-G with associated Cartan matrix $(A,\tau)$ (here we consider any Cartan matrix).  Let $\h=< h_{i}>_{i \in I}$ be the Cartan subalgebra of $\g$.
Let $\nilp_{+} \: (\text{resp., } \nilp_{-}) $ be the nilpotent Lie sub-superalgebra of $\g$ generated by $e_{i} \text{'s }$ $(\text{resp., } f_{i}\text{'s})$.  Let $\borel_{\pm} := \nilp_{\pm}\oplus \h $ be the Borel Lie sub-superalgebra of $\g$.  

Let $\eta_{\pm}:\borel_{\pm} \rightarrow \g \oplus \h$ be defined by
$$\eta_{\pm}(x)=x\oplus (\pm \bar x),$$
where $\bar x$ is the image of $x$ in $\h$.  Using this embedding we can regard $\borel_{+} $ and $\borel_{-}$ as Lie sub-superalgebras of $\g \oplus \h$

%  Def of (,)

From Proposition 2.5.3 and 2.5.5 of \cite{K} there exists a unique (up to constant factor) non-degenerate supersymmetric invariant bilinear form $(,)$ on $\g$.  Moreover, the restriction of this form to the Cartan sub-superalgebra $\h$ is non-degenerate.
 Let $(,)_{\g \oplus \h}:=(,)-(,)_{\h}$, where $(,)_{\h}$ is the restriction of $(,)$ to $\h$.

%  Proposition    $(\gl \oplus \h, \borel_{+}, \borel_{-})$ is a super Manin triple 

\begin{prop}\label{P:glManinT}
$(\g \oplus \h, \borel_{+}, \borel_{-})$ is a super Manin triple with $(,)_{\g \oplus \h}$.  
\end{prop}

\begin{pf}
Under the embedding $\eta_{\pm}$ the Lie subsuperalgebra $\borel_{\pm}$ is isotropic with respect to $(,)_{\g \oplus \h}$.   Since $(,) $ and $(,)_{\h}$ both are invariant super-symmetric nondegenerate bilinear forms then so is $(,)_{\g \oplus \h}$.  Therefore the Proposition follows.
\end{pf}

The Proposition implies that $\g \oplus \h, \borel_{+}$ and $ \borel_{-}$ are Lie superbialgebras.  Moreover, we have that $\borel_{+}^{*}\cong\borel_{-}^{op}$ as Lie superbialgebras, where $^{op}$ is the opposite cobracket.  

A straightforward calculation from the definition (see 4.4.2 of \cite{ES}) shows that $0\oplus\h$ is an ideal of the Lie superbialgebra $\g \oplus \h$.  Therefore, $\g \oplus \h / (0\oplus\h) \cong \g$ is a Lie superbialgebra.  Now from Proposition 8 of \cite{A}, we have that ($\g,\bar{r})$ is a quasitriangular Lie superbialgebra where $\bar{r}$ is the image of the canonical element $r$ in $D(\borel_{+})\cong \g \oplus \h$ under the natural projection (for the definitions of $r$ and $D$ see \ref{S:Lbialg}). 

\subsection{The Drinfeld-Jimbo type superalgebra $\DJg$} \label{SS:GenRel}  
As mentioned above, Yamane defined a QUE superalgebra for any Cartan matrix associated to the superalgebras of type A-G.  In this subsection we will summarize his results for Cartan matrices coming from root systems with exactly one odd root.

Set  $$\begin{bmatrix}
      m+n \\
      n 
\end{bmatrix}_{t}=\prod_{i=0}^{n-1}((t^{m+n-i}-t^{-m-n+i})/(t^{i+1}-t^{-i-1}))\in \C[t].$$
Let $\g$ be a  Lie superalgebra of type A-G with associated Cartan matrix $(A,\tau)$ where $\tau=\{m\}$ or $\tau=\emptyset$.  The matrix $A$ is symmetrizable, i.e. there exists nonzero rational numbers $d_{1},\dots,d_{s}$ such that $d_{i}a_{ij}=d_{j}a_{ji}$.  By rescaling, if necessary, we may and will assume that $d_{1}=1$.

%def of \DJg

Let $h$ be an indeterminate.  Set $q=e^{h/2}$ and $q_{i}=q^{d_{i}}$.  

\begin{defn}[\cite{Yam91,Yam94}] Let $\DJg$ be the $\C[[h]]$-superalgebra generated by the elements $h_{i},e_{i}$ and $ f_{i}, $ $i \in I$ satisfy the relations:
\begin{align}
\label{E:DJglRelation1}
 [h_{i},h_{j}]  &=0, & [h_{i},e_{j}]=&a_{ij}e_{j}, & [h_{i},f_{j}]=&-a_{ij}f_{j},  
\end{align}
\begin{align}
\label{E:DJglRelation2}
 [e_{i},f_{j}]=&\delta_{i,j}\frac{q^{d_{i}h_{i}}-q^{-d_{i}h_{i}}}{q_{i}-q_{i}^{-1}},  
 \end{align}
\begin{align}
\label{E:QserreA}
    e_{i}^{2}=&0   \quad  \text{for} \quad i\in I \text{ such that } a_{ii}=0,  
\end{align}
\begin{align}
\label{E:QserreA2}
    [e_{i}, e_{j}]=&0   \quad  \text{for} \quad i, j \in I \text{ such that } a_{ij}=0 \text{ and } i \neq j,  
\end{align}
\begin{align}
\label{E:QserreB}
   \sum_{v=0}^{1+ | a_{ij} |} (-1)^{v}
\begin{bmatrix}
      1 + |a_{ij}|  \\
      v  
\end{bmatrix}_{q_{i}} e_{i}^{   1 + |a_{ij}| -v}e_{j}e_{i}^{v}=0  \quad  \text{for} \quad 1 \leq i \neq j \leq s \text{ and }  i \notin \tau,
\end{align}
\begin{multline}
\label{E:QserreC}
e_{m}e_{m-1}e_{m}e_{m+1}+e_{m}e_{m+1}e_{m}e_{m-1}+e_{m-1}e_{m}e_{m+1}e_{m}
+e_{m+1}e_{m}e_{m-1}e_{m}\\-(q+q^{-1})e_{m}e_{m-1}e_{m+1}e_{m}=0  \quad \text{ if $ m-1,m, m+1\in I $ and $  a_{mm}=0$,}
\end{multline}
\begin{multline}
\label{E:QserreD}
e_{m-1}e_{m}^{3}- (q+q^{-1}-1)e_{m}e_{m-1}e_{m}^{2}
-(q+q^{-1}-1)e_{m}^{2}e_{m-1}e_{m} + e_{m}^{3}e_{m-1}=0\\ \quad \text{ if the Cartan Matrix $A$ is of type B, $\tau = \{m\}$ and $s=m$}.
\end{multline}
and the relations (\ref{E:QserreA})-(\ref{E:QserreD}) with $e$ replaced by $f$.
All generators are even except for $e_{i}$ and $f_{i}$ ($i \in \tau$) which are odd. 
\end{defn}
 We call the relations (\ref{E:QserreA})-(\ref{E:QserreD}) the \emph{quantum Serre-type relations}. 

Khoroshkin and Tolstoy \cite{KTol} and Yamane \cite{Yam91,Yam94} used the quantum double notion (see \S \ref{SS:DualDouble}) to give $\DJg$ an explicit structure of a quasitriangular Hopf superalgebra.  In the remainder of this subsection we recall some of their results which are needed in this paper.  Let $\DJb$ be the Hopf sub-superalgebra of $\DJg$ generated $\h$ and elements $e_{i}, \: i=1,\dots ,n+m-1$.  By construction $\DJb$ is a QUE superalgebra.  From Proposition \ref{P:DoubleQUE} we have that $(D(\DJb),\hat R)$ is a quasitriangular QUE superalgebra, where $\hat R$ is the canonical element of $D(\DJb)=\DJb \otimes \DJb\qd$.  There exists a epimorphism from $D(\DJb)$ to $\DJg$ coming from the identification of $\h \in \DJb$ and $\h^{*}\in \DJb\qd$.  Let $R$ be the image of $\hat R$ under this epimorphism.  The following proposition is a consequence of \cite{KTol,Yam91,Yam94}.

\begin{prop} 
The pair $(\DJg, R)$ is a quasitriangular quantization of the quasitriangular Lie bialgebra $(\g,r)$ where $r$ is the image of the canonical element associated to the super Manin triple $(\g\oplus \h, \borel_{+}, \borel_{-}) $ under the projection $\g\oplus \h \rightarrow \g$ (see \S \ref{S:Lbialg}).  In particular, the coproduct and counit given by
%$$\D a = a \otimes 1 + 1\otimes a,$$  $$\D {e_{i}}= e_{i}\otimes q^{d_{i}h_{i}} + 1 \otimes e_{i},$$ $$\D {f_{i}}= f_{i}\otimes 1 + q^{-d_{i}h_{i}} \otimes f_{i},$$ $$   \epsilon(a)=\epsilon(e_{i})= \epsilon(f_{i})=0$$ for all $ a\in \h$.
\begin{align*}
\label{}
    \D {e_{i}}= & e_{i}\otimes q^{d_{i}h_{i}} + 1 \otimes e_{i}, & \D {f_{i}}  = & f_{i}\otimes 1 + q^{-d_{i}h_{i}} \otimes f_{i},  \\
   \D a = & a \otimes 1 + 1\otimes a, &  \epsilon(a) = &  \epsilon(e_{i})= \epsilon(f_{i})=0
\end{align*}
all $ a\in \h$.
\end{prop}
\begin{pf}
From \cite{KTol,Yam91,Yam94} we have that $(\DJg, R)$ is a quasitriangular QUE superalgebra.  We only need to show that relation (\ref{E:RmatModh2}) holds.  But this follows from the explicit formula for the R-Matrix given in Theorem 8.1, and equations (8.4) and (8.5) of \cite{KTol}.  Also, see \cite{Yam94}.
\end{pf}
We call $\DJg$ the Drinfeld-Jimbo type quantization of $\g$.

%  Section      super KZ associator

% marker
\section{The superalgebra $\kzt$ and the Drinfeld category}\label{S:KZ}

In this section we define a quasitriangular quasi-Hopf superalgebra structure on $U(\g)[[h]]$.  This construction is due to Drinfeld.   We also define the Drinfeld category associated to a Lie superalgebra.  

\subsection{The quasitriangular quasi-superbialgebra $\U[[h]]$}\label{SS:BQBAU}
A superalgebra $A$ is a \emph{quasitriangular quasi-superbialgebra} if there exist even algebra homomorphisms $\Delta : A \rightarrow A\otimes A$ and $\epsilon : A \rightarrow k$ and invertible homogeneous elements $R\in A^{\otimes 2}$ and $\Phi \in A^{\otimes 3}$ such that 
\begin{equation}
\label{R:U1}
  (1\otimes \Delta)\otimes \Delta  =\Phi(\Delta \otimes 1)\otimes \Delta \Phi^{-1},
\end{equation}
\begin{equation}
\label{R:U2}
  \Delta^{op}R =R\Delta,
\end{equation}
\begin{equation}
\label{R:U3}
(1 \otimes \epsilon \otimes 1)(\Phi)=1 \otimes 1,
\end{equation} 
\begin{equation}
\label{R:Pent}
\Phi_{1,2,34}\Phi_{12,3,4}=\Phi_{2,3,4}\Phi_{1,23,4}\Phi_{1,2,3},
\end{equation}
and the hexagon relations
\begin{align}
\label{R:U4}
(\Delta \otimes 1)(R) &=\Phi_{312}R_{13}\Phi^{-1}_{132}R_{23}\Phi, & (1 \otimes \Delta)(R) &=\Phi_{231}^{-1}R_{13}\Phi_{213}R_{12}\Phi^{-1}.
\end{align}
Relation (\ref{R:U3}) implies that $\Phi$ is even.  Also, from relations (\ref{R:U4}) it follows that $(\epsilon \otimes 1)R = (1 \otimes \epsilon )R=1$ which implies that $R$ is even.

Recall the definition of the $\flip$ given in (\ref{E:Flip}). 
  Let $\g$ be a finite dimensional Lie superalgebra and let $\U$ be its universal enveloping superalgebra.  
Let $t$ be an even invariant super-symmetric element of $\g \otimes \g$, i.e. an element $t=\sum_k g_k \otimes h_k$ such that $\p g _k = \p h_k$ for all $k$,
$$\s(t)=t \text{  and  } [g\otimes 1 + 1 \otimes g,t]=0 \text{ for all } g \in \g.$$
For $n\in \N $ define $t_{ij} \in \U^{\otimes n}$ for all  $i<j$ (resp. $i>j$) by $t$ (resp. $\s(t)$) acting on the $i^{th}$ and $j^{th}$ components of the tensor product $\U^{\otimes n}[[h]]$. 

Consider the system of differential equations 
\begin{equation}
\label{E:KZ3}
\frac1\hbar \frac{\partial w}{\partial z} = \Big(\frac{t_{21}}{z} + \frac{t_{23}}{z-1}\Big)w
\end{equation}
where $\hbar=h/(2 \pi \sqrt{-1})$.
This system of equations has singularities at $0,1$ and $\infty$.  It follows from the theory of differential equations that a analytic solution on $(0,1)$ with a given initial value is unique.  Let $F_{0}(z)$ and $F_{1}(z)$ be the solutions of \eqref{E:KZ3} define on $(0,1)$ which have the asymptotic behavior $F_{0}(z) \sim z^{\hbar t_{12}}$ as $z\rightarrow 0$ and $F_{1}(z) \sim (1-z)^{\hbar t_{23}}$ as $z\rightarrow 1$.

Define $\Phi$ to be the invertible element such that $F_0(z)=F_1(z)\Phi$.  
We call $\Phi$ the \emph{super KZ associator}. 

%%
%%  Theorem   \U[[h]]   quasitriangular quasi-superbialgebra.
%%
\begin{thm}\cite{D3} \label{T:AisQuasiAlg}
The superalgebra $\big(\U[[h]], \Delta, \epsilon, \Phi, R:=e^{ht/2} \big)$ is a quasitriangular quasi-superbialgebra.
\end{thm}

\begin{pf}
In \cite{D3} Drinfeld defines a Lie algebra $\mathfrak{a_{n}}$ as the free Lie algebra with generators $X_{ij}$, $1\leq i \neq j\leq n$, module the relations 
$$ X_{ij}-X_{ji}=0,$$
$$ [X_{ij}, X_{kl}]=0,$$
$$ [X_{ij}+X_{ik}, X_{jk}]=0,$$
for $i\neq j\neq k \neq l$.  Replacing $ht_{ij}$ in the KZ-equation by $X_{ij}$, Drinfeld showed that relations  \eqref{R:U1}-\eqref{R:U4} hold.  Now let $\g$ be a Lie algebra with an $\g$-invariant symmetric two tensor $t$ in any symmetric linear tensor category.  The morphism
 $$U(\mathfrak{a_{n}}) \rightarrow U(\g^{\otimes n})$$
 given by $X_{ij}\mapsto t_{ij}$, imposes relations analogous to relations \eqref{R:U1}-\eqref{R:U4} on $U(\g)$.  Thus, applying the above discussion to the symmetric linear tensor category of superspaces the result follows.  
\end{pf}

The quasi-superbialgebra $(\U[[h]],\Delta,\epsilon,\Phi)$ is a deformation of the quasi-Hopf superalgebra $\U$, i.e. $\U[[h]]/h \U[[h]]$ is isomorphic, as a quasi-Hopf superalgebra, to $\U$.  As in the non-super case deformations of quasi-Hopf superalgebra are quasi-Hopf superalgebras (see \cite{ES}).  Therefore,  there exists a homomorphism $S: \U[[h]] \rightarrow \U[[h]]$ such that $(\U[[h]],\Delta,\epsilon,S,\Phi)$ is a quasi-Hopf superalgebra.   In summary, we have constructed a topologically free quasitriangular quasi-Hopf superalgebra $ \big(\U[[h]], \Delta, \epsilon, \Phi, R, S \big)$ which we denote by $\kzt$.

%   Drinfeld category

\subsection{The Drinfeld category} \label{SS:DriCat}

Let $\g_{+}$ be a finite dimensional Lie superbialgebra over $k$ and let $\g=D\g_{+}=\g_{+}\oplus \g_{-}$ be the Drinfeld double of $\g_{+}$ (see \S \ref{S:Lbialg}). Let $\cas$ be the Casimir element defined in (\ref{D:Cas}).  As noted $\cas$ is an even, invariant, super-symmetric element of $\g \otimes \g$.  Let $\Phi$ and $R=e^{ht/2}$ be the element arising form the pair $(\g,t)$, where $t=\cas$ (see \S\ref{SS:BQBAU}).  

Let $\M$ be the category whose objects are $\g$-modules and whose morphisms are given by $\Hom_{\M}(V,W)=\Hom_{\g}(V,W)[[h]]$.  For any $V, W \in \M$, let $V \otimes W$ be the usual super tensor product.  %Let $\beta_{V,W} : V\otimes W \rightarrow W\otimes V$ be the braiding given by the morphism $\beta = \s e^{\frac{h \cas}{2}}$ (where $\s$ is the morphism defined in subsection \ref{S:Lbialg}). 
Let $\beta_{V,W} : V\otimes W \rightarrow W\otimes V$ be the morphism given by the action of $e^{{h \cas }/2}$ on $V \otimes W$ composed with the morphism $\flip_{V,W}$ which is defined in  (\ref{E:Flip}). 
 For $V,W,U\in \M$, let $\Phi_{V,W,U}$ be the morphism defined by the action of $\Phi$ on $V \otimes W\otimes U$ regarded as an element of $\Hom_{\M}((V \otimes W)\otimes U,V \otimes (W\otimes U))$.  The morphisms $\Phi_{V,W,U}$ and $\beta_{V,W}$ define a braided tensor structure on the category $\M$ (see \cite[Prop. XIII.1.4]{Kas}), which we call the \emph{Drinfeld category}.

% 
%
%
%
%  Is this cat equivalent to A_{g,t}-mod_{fr}?
%
%
%
%
%

%section the first quantization of Lie superbialg 

% marker
\section{The quantization of Lie superbialgebras, Part I}\label{S:NQof}
In this section we give the first of two quantizations of Lie superbialgebras.  The quantization of this section is important because it commutes with taking the double.  The second quantization given in \S\ref{S:2QUofg} is important because it is functorial.  In \S\ref{S:2QUofg} we show that these two quantizations are isomorphic.  

The outline of this section is as follows.  Let $\g_{+}$ be a finite dimensional Lie superbialgebra and $\g$ its double (see \ref{S:Lbialg}).  We use Verma modules $M_{\pm}$ over $\g$ to define a forgetful functor $\OF$ from the Drinfeld category $\M$ (see \S\ref{SS:DriCat}) to the category of topologically free $k[[h]]$-modules.  We show that the endomorphisms of $\OF$ are isomorphic to a quasitriangular quantization $\OH$ of $\g$.   We then construct a Hopf sub-superalgebra $\OUh(\g_{+})$ of $\OH$ that is a quantization of $\g_{+}$.  We conclude the section with the important result that the quantum double of $\OUh(\g_{+})$ is isomorphic to $\OH$.  The last result is the main step in showing the quantization commutes with taking the double.  The results given in this section are straight forward generalizations of \cite{EK1}.

Throughout this section we use the notation of \S \ref{SS:DriCat}.
When defining maps from topologically free $\U[[h]]$-modules it is helpful to use the following isomorphism, 
\begin{equation}
\label{E:Uhmodules}
\Hom_{\U[[h]]}(X[[h]],Y)\cong \Hom_{\U}(X,Y)
\end{equation}
for any $\U$-module $X$ and topologically free $\U[[h]]$-module $Y$.

%   subsection       1st  Forgetful Functor

\subsection{The Tensor Functor $\OF$}\label{SS:TonF}
Let $M_{+}, M_{-} \in \M$ be the induced Verma modules given by
\begin{align}
\label{E:VermaMod}
   M_{+}&=\U\otimes_{_{U(\g_{+})}} c_{+} &  M_{-}=&\U\otimes_{_{U(\g_{-})}} c_{-}
\end{align}
%$$M_{+}=\U\otimes_{_{U(\g_{+})}} c_{+}$$
%$$M_{-}=\U\otimes_{_{U(\g_{-})}} c_{-}$$
where $c_{\pm}$ is the 1-dimensional trivial $\g_{\pm}$-module.  The Poincare-Birkhoff-Witt Theorem implies that the linear homomorphisms $U(\g_{+})\otimes U(\g_{-}) \rightarrow \U$ and $U(\g_{-})\otimes U(\g_{+}) \rightarrow \U$ are isomorphisms.  These isomorphisms imply that 
$$M_{\pm}=U(\g_{\mp}){1_{\pm}} $$
where ${1_{\pm}}\in M_{\pm}$  %are $\g_{\pm}$-invariant elements. 
 In particular, $M_{\pm} $ is a free $U(\g_{\mp})$-module.  
 \begin{lem}\label{L:UisoMpm}
The designation $1\mapsto 1_{+}\otimes 1_{-}$ extends to a linear map $ \phi : U(\g) \rightarrow M_{+}\otimes M_{-}$ which is an even isomorphism of $\g$-modules.
\end{lem}
\begin{pf}
By the universal property of $\U$ the linear map 
$$\g\rightarrow M_{+}\otimes M_{-}  \text{ given by } 1 \mapsto x1_{+}\otimes 1_{-}+ 1_{+}\otimes x1_{-}$$
extends to a $\g$-module morphism $\phi: U(\g) \rightarrow M_{+}\otimes M_{-}$.  By definition this morphism is even.  Moreover, it is easy to check (using the standard grading of universal enveloping superalgebras) that $\phi$ is an isomorphism.
\end{pf}
Define the functor $\OF:\M \rightarrow \A$ as
$$\OF(V)=\Hom_{\M}(M_{+}\otimes M_{-}, V)$$
where $\A$ is the category of topologically free $k[[h]]$-modules (see \S \ref{SS:DualDouble}).  As stated in \S \ref{S:Lbialg} 
the set of morphisms between superspaces is a superspace, we give $\OF(V)$ this superspace structure.    The isomorphism $\phi$ of Lemma \ref{L:UisoMpm} implies that the map 
\begin{equation}
\label{E:FVisoV[[h]]}
 \Fiso_{V} : \OF(V) \rightarrow V[[h]] \; \text{ given by } \; f \mapsto f(1_{+}\otimes 1_{-})
\end{equation}
 is a even isomorphism of superspaces.  

We now show that the functor $\OF$ is a tensor functor, i.e. there exists a family of isomorphisms $(\OJ_{V,W})_{V,W\in \M}$ such that
\begin{equation}
\label{E:Jsat}
\OJ_{U\otimes V, W}\circ (\OJ_{U,V}\otimes 1)=\OJ_{U,V\otimes W}\circ (1 \otimes \OJ_{V,W})
\end{equation}
for all $U,V,W \in \M$.  
Let $i_{\pm}: M_{\pm} \rightarrow M_{\pm} \otimes M_{\pm}$ be the ``coproduct'' on $M_{\pm}$ determined by  $i_{\pm}({1_{\pm}})={1_{\pm}} \otimes {1_{\pm}}$.  As in \cite[Lemma 2.3]{EK1} the $\g$-module morphism $i_{\pm}$ is coassociative, i.e. $(i_{\pm}\otimes 1) \circ i_{\pm}=(1 \otimes i_{\pm}) \circ i_{\pm}$ in $\Hom_{\M}(M_{\pm},M_{\pm}^{\otimes 3})$.

%   1st  DEF of \OJ_{V,W}

\begin{defn}\label{D:DefOfJ}
For each pair $V,W \in \M$ define $\OJ_{V,W}: \OF(V)\otimes \OF(W) \rightarrow \OF(V\otimes W)$ by 
$$\OJ_{V,W}(v\otimes w)=(v\otimes w) \circ \Phi^{-1}_{1,2,34} \circ (1\otimes \Phi_{2,3,4}) \circ  \beta_{23} \circ (1\otimes \Phi^{-1}_{2,3,4}) \circ \Phi_{1,2,34} \circ ( i_{+}\otimes i_{-}) .$$
\end{defn}

\begin{thm}\label{T:OJ}
The functor $\OF$ with the family $(\OJ_{V,W})_{V,W\in \M}$ is a tensor functor.
\end{thm}
\begin{pf}
Proposition 19.1 in \cite{ES} is the analogous statement in the case of Lie bialgebras.  The proof in \cite{ES} is pictorial.  It relies on the pictorial representation of $i_{\pm}$ being coassociative.  The same pictorial representation of the coassociativity holds in our case.  The proof follows exactly as in \cite{ES} after reinterpreting the pictures in our case.
\end{pf}

%    subsection   1st    The quantization of the double $\g$

\subsection{The quantization of $\g=D(\g_{+})$}\label{SS:Qofthedouble}\label{SS:EndofF}
With the use of the isomorphism given in Lemma \ref{L:UisoMpm}, 
the functor $\OF$ can be thought of as the forgetful functor
$V\mapsto \Hom_{\M}(U(\g), V).$ 
The general philosophy of tensor categories says that every forgetful functor, which is a tensor functor, induces a bialgebra structure on the underlying algebra (see \S 18.2.3 of \cite{ES}).  In this subsection, we will follow this philosophy and show that the tensor functor $\OF$ induces a superbialgebra structure on $\U[[h]]$.  We do this in three steps: (1) show that endomorphisms of $\OF$ are isomorphic to $U(\g)[[h]]$, (2) show that the family $(\OJ_{V,W})_{V,W\in \M}$ is determined by an element $\OJ\in U(\g)[[h]]^{\otimes 2} $, (3) use $\OJ$ to define a quasitriangular Hopf superalgebra structure on $U(\g)[[h]]$.

Let $\End(\OF)$ be the algebra of natural endomorphisms of $\OF$.  In other words, $\End(\OF)$ is the algebra consisting of elements $\eta$, so that each $\eta$ is a collection of linear morphisms
$\eta_{V} : \OF(V) \rightarrow \OF(V)$
such that for all $V,W \in \M$ and $f : V \rightarrow W$ we have $\OF(f) \circ \eta_{V}=\eta_{V} \circ \OF(f)$.  We make $\End(\OF)$ a superspace by defining $\eta \in(\End(\OF))_{\p i}$ if the parity of $\eta_{V}$ is $\p i$ for all $V\in \M$.  This makes $\End(\OF)$ into a superalgebra.

\begin{lem}\label{L:UHF}
	There is a canonical even superalgebra isomorphism
	\begin{equation}
\label{E:U[[h]]isoEndF}
\theta : \U[[h]] \rightarrow  \End(\OF) \text{ given by } x \mapsto x_{|V}
\end{equation}
	where $x_{|V}$ is $x$ acting on the $\U[[h]]$-module $V[[h]]$.
\end{lem}

\begin{pf}
	Using the even isomorphism (\ref{E:FVisoV[[h]]}) we can identify $\OF(V)$ and $V[[h]]$.  Under this identification $\theta (x)=x_{|V}\in \End(\OF)$ is the endomorphism given by the action of $x$ on $V[[h]]$.  The morphism $\theta$ is even since the action of a homogeneous element $x$ on $V[[h]]$ preserves the grading.  If $x\neq y$ then $x_{|\U}\neq y_{|\U}$ implying $\theta$ is one to one.  
	
	Next we will show that $\theta$ is onto.  Let $\eta \in \End(\OF)$, using the above isomorphism we will think of $\eta_{V}$ as a map from $V[[h]]$ to itself.  Set $x = \eta_{\U}(1)$.  Let $y\in \U$ and let $r_{y}$ be the element of $ \End(\U)$ given by 
	$$ r_{y}(z)=(-1)^{\p y \p z}zy$$
for $z \in \U$.  %In other words, $r_{y}$ is right multiplication by $y$ with the appropriate sign.  
 Note that $\OF(r_{y})$ under the isomorphism $\OF(\U) \rightarrow \U[[h]]$ is $r_{y}$.  We have
$$\eta_{\U}(y)=\eta_{\U}(r_{y}1)=(-1)^{\p y \p x}r_{y}\eta_{\U}(1)=(-1)^{\p y \p x}r_{y}x=xy.$$
Combining this calculation with (\ref{E:Uhmodules}), we have $\eta_{\U}=l_{x}$ where $l_{x}(z)=xz$ for $z \in \U$.  
Similarly $\eta_{V}=x_{|V}$ for any free $\g$-module $V$.  This shows that $\theta$ is onto since every $\g$-module is a quotient of a free module.
\end{pf}

%From now on we will make no distinction between $\U[[h]]$ and $\End(\OF)$, identifying them by $\theta$.

In the rest of this subsection we use properties of the tensor functor $\OF$ and the isomorphism $\theta$ to put algebraic structures on $\U[[h]]$.  
 
Define the element $\OJ \in \U^{\otimes 2}[[h]]$ to be
\begin{equation}
\label{E:defJ}
\OJ=(\phi^{-1} \otimes \phi^{-1})\left(  \Phi^{-1}_{1,2,34}  (1\otimes \Phi_{2,3,4})    \beta_{23}  (1\otimes \Phi^{-1}_{2,3,4})  \Phi_{1,2,34}  \left( 1_{+}\otimes 1_{+} \otimes 1_{-} \otimes 1_{-} \right)\right)
\end{equation}
where $\phi$ is the isomorphism given in Lemma \ref{L:UisoMpm}.   %The following lemma shows that $\theta(\OJ)$ and $\OJ_{V,W}$ are essentially the same after the identification of $\OF(V)$ with $V[[h]]$. 

\begin{lem}\label{L:JandJvw}
Let $\theta$ be the isomorphism of Lemma \ref{L:UHF}.  Then $\theta(\OJ)=\OJ_{V,W}$, i.e. 
\begin{equation}
\label{E:JandJvw}
\OJ(v\otimes w)=\Fiso_{V\otimes W}(\OJ_{V,W}(\Fiso_{V}^{-1}(v)\otimes \Fiso_{W}^{-1}(w)))
\end{equation}
for $v\in V[[h]]$ and $w\in W[[h]]$
\end{lem}
\begin{pf}
For each $v\in V[[h]]$ let $f_{v}$ to be the element of $\OF(V)$ defined by $f_{v}(x)=(-1)^{\p v \p x}v$ for $x \in M_{+}\otimes M_{-}$.  Notice that the element $f_{v}$ has parity $\p v$.  From Lemma \ref{L:UisoMpm} we have $f_{v}(1_{+}\otimes 1_{-})=v$ which implies that $f_{v}=\Fiso_{V}^{-1}(v)$.  To simplify notation let
$$\Otheta_{1}\otimes \Otheta_{2} =  \Phi^{-1}_{1,2,34}  (1\otimes \Phi_{2,3,4})    \beta_{23}  (1\otimes \Phi^{-1}_{2,3,4})  \Phi_{1,2,34}  \left( 1_{+}\otimes 1_{+} \otimes 1_{-} \otimes 1_{-} \right)$$
be the element of $ (M_{+}\otimes M_{-})^{\otimes 2}[[h]]$.
Now we have the right side of (\ref{E:JandJvw}) is
\begin{align*}
\label{}
  \left(\OJ_{V,W}(\Fiso_{V}^{-1}(v)\otimes \Fiso_{W}^{-1}(w))\right)(1_{+}\otimes 1_{-})  
    & =  (f_{v}\otimes f_{w})(\Otheta_{1}\otimes \Otheta_{2})\\
    & = (-1)^{\p w \p{\Otheta}_{1}}f_{v}\Otheta_{1}\otimes f_{w}\Otheta_{2}\\
    & = (-1)^{\p w \p{\Otheta}_{1}+\p v \p{\Otheta}_{1}+ \p w \p{\Otheta}_{2}}\phi^{-1}(\Otheta_{1})v \otimes \phi^{-1}(\Otheta_{2})w.
\end{align*}
On the other hand, the left side of (\ref{E:JandJvw}) is
\begin{align*}
(\phi^{-1}(\Otheta_{1})\otimes \phi^{-1}(\Otheta_{2}))(v\otimes w)
	& = (-1)^{\p{\Otheta}_{2} \p v}\phi^{-1}(\Otheta_{1})v\otimes \phi^{-1}(\Otheta_{2})w\\
	& = (-1)^{\p w \p{\Otheta}_{1}+\p v \p{\Otheta}_{1}+ \p w \p{\Otheta}_{2}}\phi^{-1}(\Otheta_{1})v \otimes \phi^{-1}(\Otheta_{2})w.
\end{align*}
where the last equality follows from the fact that $\Otheta_{1}\otimes \Otheta_{2}$ is even, i.e. $\p \Otheta_{1}+\p \Otheta_{2}=0$.
\end{pf}

%     1st   Lemma      \OJ = 1 + \frac{rh}{2} \mod h^{2}

\begin{lem} \label{L:Jmod1}
$\OJ \equiv 1 + \frac{rh}{2} \mod h^{2}$.
\end{lem}
\begin{pf}
Recall the definition of $r$ given in \S\ref{S:Lbialg}, i.e. $r=\sum_{i}p_{i}\otimes m_{i}$ where $(p_{i})_{i}$ and $(m_{i})_{i}$ are a bases of $\g_{+}$ and $\g_{-}$ respectively.  It follows that 
\begin{equation}
\label{E:rflipaction}
\flip(r)(1_{-}\otimes 1_{+})=\sum (-1)^{\p m_{i}}m_{i}\otimes p_{i})(1_{-}\otimes 1_{+})=0,
\end{equation}
\begin{equation}
\label{E:actionofpi}
p_{i}\phi^{-1}(1_{+}\otimes 1_{-})= \phi^{-1} ( (p_{i}1_{+})\otimes 1_{-} + 1_{+}\otimes (p_{i} 1_{-}))= \phi^{-1}(1_{+}\otimes (p_{i} 1_{-}))
\end{equation}
where $\phi$ is the isomorphism given in Lemma \ref{L:UisoMpm}.

From the hexagon relation (\ref{R:U4}) we have $\Phi \equiv 1 \mod h^{2}$.  Thus,
\begin{align*}
\label{}
   \OJ & % \equiv  (\phi^{-1} \otimes \phi^{-1} )(1\otimes \flip \otimes 1)(e^{h\cas_{23}/ 2})(1_{+}\otimes 1_{+}\otimes 1_{-}\otimes 1_{-}) && \mod h^{2} \\
     \equiv  (\phi^{-1}  \otimes \phi^{-1} )(e^{h\cas_{23}/ 2})(1_{+}\otimes 1_{-}\otimes 1_{+}\otimes 1_{-}) && \mod h^{2}  \\
    & \equiv  1 + h/2  (\phi^{-1}  \otimes \phi^{-1} )(r_{23} +\flip (r)_{23} )(1_{+}\otimes 1_{-}\otimes 1_{+}\otimes 1_{-}) &&  \mod h^{2}  \\
    & \equiv  1 +h/2  (\phi^{-1}  \otimes \phi^{-1} )(\sum 1_{+}\otimes p_{i}1_{-}\otimes m_{i}1_{+}\otimes 1_{-}) && \mod h^{2}  \\
    & \equiv   1 + h/2 \sum (p_{i}\phi^{-1}(1_{+}\otimes 1_{-})  \otimes m_{i}\phi^{-1} ( 1_{+}\otimes 1_{-}) ) && \mod h^{2}  \\
    & \equiv 1 + rh/2 && \mod h^{2}  
\end{align*}
where the third equivalence follows from (\ref{E:rflipaction}) and the fourth follows from (\ref{E:actionofpi}).
\end{pf}
	
\begin{prop}\label{P:OHisTwistJ}
Let $\OH=\U[[h]]$.  Then $\OH$ is a Hopf superalgebra whose coproduct, counit and antipode are given by
\begin{equation}\label{E:DefCoproduct}
\begin{array}{cc}
   \Delta = \OJ^{-1}\Delta_{0}\OJ,        &    \epsilon=\epsilon_{0}
\end{array}
\end{equation}
\begin{equation}\label{E:DefAntip}
S=QS_{0}Q^{-1}
\end{equation}
where $Q=m(S_{0}\otimes 1 )(\OJ)$ and $\Delta_{0}$, $\epsilon_{0}$ and $S_{0}$ are the usual coproduct, counit and antipode of $\U[[h]]$.  
\end{prop}	
\begin{pf}
First, $\Delta$ and $\epsilon$ are algebra morphisms since $\Delta_{0}$ and $\epsilon_{0}$ are algebra morphisms. 
From Lemma \ref{L:Jmod1} we have that  $(\epsilon \otimes 1)\OJ = (1 \otimes \epsilon)\OJ = 1$ which implies $(\epsilon \otimes 1)\Delta =1=(1 \otimes \epsilon)\Delta$.  Theorem \ref{T:OJ} and Lemma \ref{L:JandJvw} imply that 
\begin{equation}
\label{E:OJ123}
\OJ_{12,3}(\OJ \otimes 1) = \OJ_{1,23}(1\otimes \OJ). 
\end{equation} 
We will now use equality (\ref{E:OJ123}) to show that $\Delta$ is coassociative.  
\begin{align}
\label{}
  (1 \otimes \Delta) \Delta(x) 
   & = (1 \otimes \OJ^{-1}\Delta_{0}\OJ)(\OJ^{-1}\Delta_{0}(x)\OJ)\\
   & = (1 \otimes \OJ^{-1})\OJ_{1,23}^{-1}(1\otimes \Delta_{0})\Delta_{0}(x)\OJ_{1,23}(1 \otimes \OJ)\\
   & = (\OJ^{-1} \otimes 1)\OJ_{12,3}^{-1}(\Delta_{0} \otimes 1)\Delta_{0}(x)\OJ_{12,3}(\OJ\otimes 1)\\
   & = ( \Delta \otimes 1 ) \Delta(x)
\end{align}
for all $x\in \OH$. The compatibility conditions between $S$ and $\epsilon$ follow in a similar manner.  
\end{pf}

The isomorphism $\theta$ of Lemma \ref{L:UHF} induces a Hopf superalgebra on $\End(\OF)$.  For the  rest of this paper, we identify the Hopf superalgebra $\OH$ with $\End(\OF)$ (using $\theta$).  As we will see it is sometimes convenient to use the elements of $\OH$ and other times endomorphisms of $\End(\OF)$.  
%It is important to observe that there is a corresponding Hopf superalgebra structure on $\End(\OF)$.  For example,  let $\Delta: \End(\OF) \rightarrow \End(\OF)^{\otimes 2}$ be the linear morphism given by \begin{equation}\label{E:coprodEnd}\Delta(\eta)_{V,W}(v\otimes w)=\OJ_{V,W}^{-1}\eta_{V,W}\OJ_{V,W}(v\otimes w)\end{equation}for $\eta\in \End(\OF)$, $v\in \OF(V)$, $w\in \OF(W)$.  Theorem \ref{T:OJ} implies that $\Delta$ is coassociative.  One can continue and define a natural Hopf superalgebra on $\End(\OF)$ (see .  Under this identifcation, Lemma \ref{L:JandJvw} implies that the morphism given in (\ref{E:coprodEnd}) is the coproduct on $\OH$.  Therefore, $\Delta$ is a coproduct on $\End(\OF)$.  One can continue and define a natural  The isomorphism $\theta$ allows us to impose a Hopf superalgebra structure on $\End(\OF)$.  It should be noted that this is not the natural The isomorphism $\theta$ allows us to impose a Hopf superalgebra structure on $\End(\OF)$.  Throught the rest of this paper we identify We will use the isomorphism (\ref{E:U[[h]]isoEndF}) to identify 

%   Theorem      1st   \OH is a quantization of the super Lie bialgebra \g

\begin{thm}
\label{T:Qofg}
	$\OH$ is a quantization of the Lie superbialgebra \g.
\end{thm}
\begin{pf}
 By definition $\OH / h\OH $ is isomorphic to the Hopf superalgebra $\U$.  To prove the theorem we show that relation (\ref{R:ConditionDelta}) holds.  From the definition of the coproduct $\Delta$ and Lemma~\ref{L:Jmod1} we have
\begin{equation}
\label{E:Qofg1}
\D x \equiv \Delta_{0}(x) + (h/2) [\Delta_{0}(x),r] \mod h^{2}
\end{equation}
for all $x \in \g \subset \OH$.
Thus,
\begin{align}
\label{E:Qofg2}
  h^{-1}( \D x - \Delta^{op}(x) )& \equiv  h^{-1}\Delta_{0}(x) +   1/2 [\Delta_{0}(x),r] \nonumber \\
  & \hspace{40pt} -   h^{-1}\Delta^{op}_{0}(x) -  1/2 [\Delta^{op}_{0}(x),\flip(r)]     \mod h \nonumber \\
    &  \equiv  1/2 [\Delta_{0}(x),r-\flip_{\g,\g}(r)]  \mod h  \nonumber\\
    &  \equiv  [\Delta_{0}(x),r]    \mod h
\end{align}
since $t=r+ \flip(r)$ is \g -invariant and $\Delta_{0}(x)=\Delta_{0}^{op}(x)$ for all $x\in \g$  (for the definition of $\flip$, see (\ref{E:Flip})).  The proof is completed by recalling that the cobracket of $\g$ is defined by $\partial r(x):= [\Delta_{0}(x),r] $.
\end{pf}

%   Def of R

Define $R=(\OJ^{op})^{-1}e^{\frac{h\cas}{2}}\OJ\in \OH \otimes \OH$.  We call $R$ the R-matrix.  

%  corollary       (\OH,R)    is a quasitriangular

\begin{cor}\label{C:quasiTriQuant}
$(\OH,R)$ is a quasitriangular quantization of $(\g,r)$.
\end{cor}
\begin{pf}
Replacing the standard commutativity isomorphism with the super commutativity isomorphism, (i.e. substituting $\flip$ for $\sigma$) the proof follows just as in the purely even case \cite[Corollary 19.1]{ES}.
\end{pf}

%    subsection    1st  Quantization of $\g_{+}$ and $\g_{-}$

\subsection{Quantization of $\g_{+}$ and $\g_{-}$}\label{SS:Qofgpgm}
Here we construct a quantization of the Lie superbialgebra $\g_{\pm}$, which is a Hopf sub-superalgebra of $\OH$.  To this end, we continue following the work of Etingof and Kazhdan \cite{EK1} and notice that $R$ is polarized,  i.e. $R\in \OUhp \otimes \OUh(\g_{-}).$   It is possible to show directly that $\OUhp$ is closed under coproduct.  However, in \S \ref{SS:quantumDual} we use the polarization of $R$ to show that the quantization commutes with the double.  

%Def of U_{h}(\g_{+})

Using the even isomorphisms (\ref{E:FVisoV[[h]]}) and (\ref{E:U[[h]]isoEndF}) we can identify the superalgebras $\End(\OF)$ and $\End(M_{+} \otimes M_{-})$.  We will not make a distinction between these superalgebras.  Define $\OUhp = \OF(M_{-})$ and embed it into $\OH$ using the map $i : \OF(M_{-}) \rightarrow \End(M_{+} \otimes M_{-})$ given by
$$i(x)=(1\otimes x)\circ \Phi \circ (i_{+}\otimes 1)$$
for $x\in \OF(M_{-})$.  

The coassociativity of $i_{+}$ implies that for $x,y \in \OF(M_{-})$
$$ i(x)\circ i(y)=i(z)$$
where $z= x\circ (1 \otimes y )\circ \Phi \circ (i_{+}\otimes 1) \in \OF(M_{-})$.
Using the embedding $i$, we consider $\OUhp $ is a subsuperalgebra of $\OH$.  Similarly, the map $\OF(M_{-})\rightarrow \End(M_{+}\otimes M_{-})$ given by $x\mapsto (x\otimes 1) \circ \Phi \circ (1 \otimes i_{-})$ makes $\OUh(\g_{-}):=\OF(M_{+})$ into a subsuperalgebra of $\OH$.

%    1st U_{h}(\g_{+})  is a Hopf sub-superalgebra of $\OH$.

\begin{thm}
\label{T:UhpH}
  $\OUh(\g_{+})$ and $\OUh(\g_{-})$ are Hopf sub-superalgebra of $\OH$.  Moreover $\OUh(\g_{\pm})$ is a quantization of the Lie superbialgebra $\g_{\pm}$.
\end{thm}

\begin{pf}
As in \cite{EK1} we needed the following lemma to prove the theorem.  

%  Lemma   $R$ is polarized

\begin{lem}\label{L:Rpolar}
$R$ is polarized, i.e. $R\in \OUhp \otimes \OUh(\g_{-}) \subseteq \OH \otimes \OH.$
\end{lem} 
\begin{pf}
In \cite{ES} the analogous statement in the purely even case is proved using a pictorial proof.  After representing the Hopf superalgebra structure of $\OH$ and functoriality of the braiding $\beta = \flip e^{\frac{h \cas}{2}}$ pictorially the proof follow exactly as in Lemma 19.4 \cite{ES}.  %Using the fact that $\OH$ is a Hopf superalgebra and the fact that the braiding $\beta = \s e^{\frac{h \cas}{2}}$ is functorial one can show that    pictorial relations hold. Following \cite{EK} one can interpreting the analogous pictorial r   Lemma 19.4 in \cite{ES} is the analogous statement in the purely even case.  The proof follows from the universal proof of Lemma 19.4 in \cite{ES}.
\end{pf}

%   Def of \rmap and image of rmap

Let $\rmap_{+}: \OUh(\g_{-})^{*}\rightarrow \OUh(\g_{+})$ and $\rmap_{-}: \OUh(\g_{+})^{*}\rightarrow \OUh(\g_{-})$ be the even linear maps given by
\begin{equation*}
\label{E:DefOfrmap}
 \rmap_{+}(f)= (1\otimes f)(R)  \: \text{ and } \:   \rmap_{-}(f)= (f \otimes 1)(R)  
\end{equation*}
for $f \in \OUh(\g_{\mp})^{*}:=\Hom_{\A}(\OUh(\g_{\mp}),k[[h]])$.  Let $\rmapim_{\pm}$ be the images of $\rmap_{\pm}$.  Let $\widetilde{\rmapim}_{\pm} $ be the closer of the $k[[h]]$-superalgebra generated by $\rmapim$. 

%   Lemma   the image of rmap and U_{h}(g_{\pm})

\begin{lem}
\label{L:rmapim}
$\widetilde{\rmapim}_{\pm} \otimes_{k[[h]]}k((h)) $ is the $h$-adic completion of $ \OUh(\g_{\pm}) \otimes_{k[[h]]}k((h))$ where the tensor product is the tensor product in the h-adic completion.
\end{lem}
\begin{pf}
Using the even graded linear map $\rmap_{\pm}$ the proof is identical to the proof of Proposition 4.5 \cite{EK1}.  In particular, no new signs are introduced in the proof and grading is preserved since it is preserved by $\rmap_{\pm}$.
\end{pf}
Now we prove the theorem.  Relations (\ref{R:comultR}) imply that $\rmapim_{\pm}$ is closed under coproduct.  Therefore, by Lemma \ref{L:rmapim} we have that $\OUh(\g_{\pm})$ is closed under coproduct.  Moreover, since $\OH$ is a quasitriangular Hopf superalgebra we have 
$$(S\otimes 1)R=R^{-1}$$
which implies that $\OUh(\g_{\pm})$ is closed under the antipode.  This proves $\OUh(\g_{\pm})$ is a Hopf sub-superalgebra of $\OH$.  

Next we show that $\OUh(\g_{\pm})$ is a quantization of $\g_{\pm}$.  The isomorphism given in (\ref{E:FVisoV[[h]]}) implies that $\OUh(\g_{\pm})$ is isomorphic, as a superspace, to $U(\g_{\pm})[[h]]$.  Moreover, the Hopf superalgebra $\OUh(\g_{\pm})/h\OUh(\g_{\pm})$ is isomorphic to $U(\g_{\pm})$.  Since $\OUh(\g_{\pm})$ is a Hopf sub-superalgebra of $\OH$ we have that equivalencies (\ref{E:Qofg1}) and (\ref{E:Qofg2}) hold for all $x \in \g_{\pm}$.  Thus,  $\OUh(\g_{\pm})$ satisfies equivalence (\ref{R:ConditionDelta})  and so is a quantization of $\g_{\pm}$.  
\end{pf}
We call $\OUhp$ the Etingof-Kazhdan quantization of $\g_{+}$.

%     subsection              The quantum dual of $\OUhp$

\subsection{The quantum dual of $\OUhp$}\label{SS:quantumDual}

%  We denote the quantum dual if \qd

Recall the definitions of the quantum dual and double of a QUE superalgebra given in \S\ref{SS:DualDouble}.   
In this subsection we will show that the quantum dual of $\OUh(\g_{-})^{op}$ is $\OUh(\g_{+})$ and that the double of $\OUhp$ is $\OH$.  The former statement follows from the use of the linear map $\rmap_{+}$ which arises from the polarization of $R$.  In \S \ref{S:FunctQ} we will use the results of this subsection to show that the quantization commutes with taking the double.
 
%   Theorem   $(U_{h}(\g_{-})^{op})^{*} \rightarrow U_{h}(\g_{+})$

\begin{prop}\label{P:rmapIshomo}
The linear map $\rmap_{+}$ ($\rmap_{-}$) is a even injective homomorphism of topological Hopf superalgebras $(\OUh(\g_{-})^{op})^{*} \rightarrow \OUh(\g_{+})$ (resp. $\OUh(\g_{+})^{*} \rightarrow \OUh(\g_{-})^{op}$).  Moreover,  $\rmapim_{\pm}=\OUh(\g_{\pm})^{\prime}$.
\end{prop}
\begin{pf}
The proof follows as in the proof of Proposition 4.8 and Proposition 4.11 in \cite{EK1}.
\end{pf}

\begin{cor}\label{C:qd}
The quantum dual of the QUE superalgebra $\OUh(\g_{+})$  is $\OUh(\g_{-})^{op}$.  Moreover, the  quantization of $\g=D(\g_{+})$ given in \ref{SS:Qofthedouble} and the quantum double of $\OUh(\g_{+})$ are isomorphic as quasitriangular QUE superalgebras, i.e. $\OH \cong D(\OUhp)$ (for the definition of the quantum double see Proposition \ref{P:DoubleQUE}).
\end{cor}
\begin{pf}
The first assertion follows from 
\begin{equation}
\label{E:Ug+QDual}
\OUh(\g_{+})\qd:=(\OUh(\g_{+})^{*})^{\vee} \cong ((\OUh(\g_{-})^{op})^{\prime })^{\vee}=\OUh(\g_{-})^{op}
\end{equation}
where the isomorphism comes from Proposition \ref{P:rmapIshomo} and the third equality follows from (\ref{E:QDualandDual}). %?Similarly, one shows that $(\OUh(\g_{-})^{op})\qd \cong \OUh(\g_{+})$?

To prove the second assertion we will show that $\OH$ satisfies the defining relations of the Hopf superalgebra structure on the double of $\OUh({\g_{+}})$, then the result follows from the uniqueness of Proposition \ref{P:DoubleQUE}.  
By Theorem \ref{T:UhpH}, $\OUh(\g_{+})$ and $\OUh(\g_{-})$ are Hopf sub-superalgebras of $\OH$.  The multiplication map $\OUh(\g_{+})\otimes \OUh(\g_{-}) \rightarrow \OH$ is a bijection, as it is modulo $h$.
In equation (\ref{E:Ug+QDual}) we concluded that the map $\rmap_{-}$ induces an isomorphism between $\OUh(\g_{+})\qd $ and $ \OUh(\g_{-})^{op}$.  Therefore, the definition of the quantum dual implies  
\begin{equation}
\label{E:DoubleIsEuaal}
D(\OUh(\g_{+}))\cong \OUh(\g_{+})\otimes \OUh(\g_{-}).
\end{equation}
Recall that the map was defined $\rmap_{-}$ was defined using the R-Matrix $R$ of $\OH$.  It follows that $R$ corresponds to the canonical element $\hat R$ of $D(\OUh(\g_{+}))$.  Thus, the uniqueness of Proposition \ref{P:DoubleQUE} implies $D(\OUh(\g_{+}))=\OH$.
\end{pf}

%  Quantization of quasitriangular Lie superbialgebras

% marker
\section{Quantization of quasitriangular Lie superbialgebras}\label{S:QUofquasi}
Let $\g_{+}$ be a finite dimensional Lie superbialgebra.  Recall that the double $D(\g_{+})$ of $\g_{+}$ is a quasitriangular Lie superbialgebras (see \S \ref{S:Lbialg}). 
In this section we will construct a quantization of quasitriangular Lie superbialgebras. 
This quantization is similar to the quantization of Lie superbialgebra of \S \ref{S:NQof}.  In \S \ref{S:FunctQ} we will show that for finite dimensional quasitriangular Lie superbialgebras the two quantizations are isomorphic.  Moreover, by construction the quantization $\OH$ given in \S \ref{S:NQof} is the same as the quantization of $D(\g_{+})$ given below.    These facts are used in proving that the quantization commutes with taking the double.
  
Let $(\g, r)$ be a quasitriangular Lie superbialgebra.  Set 
$$    \g_{+}=\{(1 \otimes f) r \: | \: f \in \g^{*}\}\: \text{ and } \:
   \g_{-}=\{(f \otimes 1) r \: | \: f \in \g^{*}\}.  $$
Then $\g_{+}$ and $\g_{-}$ are finite dimensional Lie superbialgebras (see \cite[Lemma 5.2]{EK1}).  Moreover, $\g_{-}\cong \g_{+}^{*}$ and there is a natural homomorphism of quasitriangular Lie superalgebras $\pi : D(\g_{+}) \rightarrow \g$ (see \S 5 of \cite{EK1}).  Let $\M$ be the category whose objects are $\g$-modules and whose morphisms are given by $\Hom_{\M}(V,W)=\Hom_{\g}(V,W)[[h]]$.  As in \S \ref{SS:DriCat}, let  $\mathcal{M_{D(\g_{+})}}$ be the Drinfeld category of the double $D(\g_{+})$.  From the homomorphism $\pi$ we have
$$\pi^{*}: \M \rightarrow  \mathcal{M_{D(\g_{+})}}$$
whose pull back gives a braided tensor structure on the category $\M$.  Let $M_{-}$ and $M_{+}$ be the ``Verma'' modules in $ \mathcal{M_{D(\g_{+})}}$ (see (\ref{E:VermaMod})).  

Let $\OFq: \M \rightarrow \A$ be the functor given by $$\OFq(V)=\Hom_{\mathcal{M_{D(\g_{+})}}}(M_{+}\otimes M_{-}, \pi^{*}(V)).$$  Then $\OFq$ is a tensor functor with the isomorphism of functors $\OJ$ giving in Definition \ref{D:DefOfJ}.  As in (\ref{E:FVisoV[[h]]}) the map 
\begin{equation}
   \OFq(V) \rightarrow V[[h]] \; \text{ given by } \; f \mapsto f(1_{+}\otimes 1_{-})
\end{equation}
is an even isomorphism of superspaces.  Using this isomorphism we construct the canonical isomorphism  $\theta : \U[[h]] \rightarrow  \End(\OFq)$ of Lemma \ref{L:UHF}.  The equations (\ref{E:DefCoproduct}) and (\ref{E:DefAntip}) define a Hopf superalgebra structure on $U(\g)[[h]]$ which is equal to $\End(\OFq)$.  Finally, as in Corollary \ref{C:quasiTriQuant} we have that $(U(\g)[[h]],R)$ is a quasitriangular quantization of $(\g,r)$, where the R-matrix is defined as in \S \ref{SS:Qofthedouble}.  We denote this quasitriangular QUE superalgebra by $\Uqt$.

%      2nd Quantization of Lie superbialgebras

% marker
\section{The quantization of Lie superbialgebras, Part II}\label{S:2QUofg}

Here we give the second quantization of Lie superbialgebras.   As mentioned before, this quantization is isomorphic to the first quantization of the Lie superbialgebra constructed in  \S\ref{S:NQof}.  We denote the quantization of this section by $\FUh(\g)$.  In \S\ref{S:FunctQ} we will see that the quantization of this section is functorial.

We follow the quantization of Lie bialgebras given in Part II of \cite{EK1}.  The results of \cite{EK1} should generalize to the setting of all Lie superbialgebras.  However, we will only check that the results hold for finite dimensional Lie superbialgebras.  %This simplifies some of the computations given in \cite{EK1}.  For example, since $\g$ is finite dimensional all $\g$-modules are equicontinuous (see \cite[\S 7.3]{EK1}).   %By construction, the quantization of Lie bialgebras given in Part II of \cite{EK1} restricted to finite dimensional Lie bialgebras is a special case of the quantization given in this section. 

In this section we consider topological superspaces.  We need topology to deal with convergence issue involving duals of infinite dimensional space and tensor products of such spaces.  In particular, we need modules to be equicontinuous (see \cite[\S 7.3]{EK1}).  Since we are working with finite dimensional Lie superbialgebras all modules are over such superalgebras are equicontinuous.  For this reason, we will assume that all modules are equicontinuous.  We proceed in much the same way as in \S \ref{S:NQof}.  In other words, given a finite dimensional Lie superbialgebra $\g_{+}$, we use Verma modules to define a tensor functor such that the set of endomorphisms of this functor is a quantization of the double of $\g_{+}$ which contains a quantization of $\g_{+}$.

\subsection{Topological superspaces}
Let $k$ be a field of characteristic zero.  We consider $k$ as a topological superspace concentrated in degree $\p 0$ with discrete topology.

Let $V$ be a topological superspace, i.e. a $\Z_{2}$-graded topological vector space.  We say $V$ is linear if open superspaces of $V$ form a basis of neighborhoods of $0$.  The superspace $V$ is separated (complete) if the natural map $V \rightarrow \underleftarrow{\lim} V / U$ is a monomorphism (resp. epimorphism) where the limit runs over open sub-superspaces $U$.  Throughout this section we will only consider complete, separated topological superspaces, so when we use the phrase ``topological superspace'' we will mean ``complete, separated, linear topological superspaces''.   

Let $V$ and $W$ be topological superspaces.  If $U$ is an open sub-superspace of $V$ then $V/U$ is discrete.   Using this we define the tensor product of two topological superspaces $V$ and $W$ to be
$$V \ttp W :=  \underleftarrow{\lim} V/ V' \otimes W/W'$$
where $V'$ and $W'$ run over open sub-superspaces of $V$ and $W$ respectively.  Let $V[[h]]=V\ttp k[[h]]$ be the space of formal poser series in $h$. 
We give the superspace $\Hom_{k}(V,W)$ of all continuous homomorphisms a topology, as follows.  Let $B$ be a topological basis of W.  For any $n \geq 1$ let $U_{1}, U_{2},...,U_{n}\in B$ and $v_{1}, v_{2},...,v_{n}\in V$.  Then the collection 
$$\left( \{ f \in \Hom_{k}(V,W) : f(v_{i})\in U_{i} \text{ for } i = 1,...,n \} \right)_{U_{1}, U_{2},...,U_{n},v_{1}, v_{2},...,v_{n}}$$
is a basis for the topology on $\Hom_{k}(V,W)$.  We call this topology the weak topology.  Note that if $V$ is finite dimensional then the weak topology on $V^{*}=\Hom_{k}(V,k)$ is the discrete topology.  

\subsection{Topological $\g$-modules} Let $\g_{+}$ be a finite dimensional Lie superbialgebra we give $\g_{+}$ the discrete topology. 
In this section a $\g_{+}$-module will be a topological superspace $M$ with a continuous homomorphism of topological Lie algebras 
$$\pi: \g_{+} \rightarrow \End(M)$$
 such that $\pi((\g_{+})_{\p i}) \subset \End(M)_{\p i}$ for $\p i =\p 0, \p 1$.
 
Let $\g= D(\g_{+})$ be the Drinfeld double of $\g_{+}$ (see \S \ref{S:Lbialg}).   Given two topological $\g$-modules $V,W$ let $\Hom_{\g}(V,W)$ be the topological superspace of all continuous $\g$-modules homomorphisms.  
Let $\Mt$ be the category whose objects are $\g$-modules and morphism are given by
$$\Hom_{\Mt}(V,W)=\Hom_{\g}(V,W)[[h]]$$
for $V,W \in \Mt$.  

Using the tensor produce $\ttp$ we define a braided tensor structure on $\Mt$ as follows.  Let $\cas$ be the Casimir element defined in \S\ref{S:Lbialg} and $\Phi$ be the associator constructed using $\cas$ (see \S \ref{SS:BQBAU}).  For $V,W,U \in \Mt$, let $\Phi_{V,W,U} $ be the element of $\Hom_{\Mt}((V \ttp W)\ttp U,V \ttp (W\ttp U))$ given by the action of $\Phi$ on $V \ttp W\ttp U$ and let $\beta_{V,W}:=\flip_{V,W} e^{{h \cas }/2} \in \Hom_{\Mt}(V \ttp W, W \ttp V)$ (for the definition of $\flip$, see (\ref{E:Flip})).  The morphisms $\Phi_{V,W,U}$ and $\beta_{V,W}$ define a braided tensor structure on $\Mt$. 

Let $\At$ be the category whose objects are $k[[h]]$-modules and morphisms are continuous $k[[h]]$-linear maps.  $\At$ is a symmetric tensor category where the tensor product $V \atp W$ is the tensor product $V \ttp W$ modulo the image of the operator $h\otimes 1-1 \otimes h$.     

Recall the definitions of $M_{\pm}$ and $i_{\pm}$ given in \S \ref{SS:TonF}.  We give $M_{\pm}$ the discrete topology.  For finite dimensional Lie bialgebras these topologies are the same as the topologies defined in \cite[\S 7.5]{EK1}.  % To define a topology on $M_{+}$ we identify $M_{+}$ with $U(\g_{-})$.  For any $n \gqe 0$ let $U(\g_{-})_{n}$ be the elements of $U(\g_{-})$ with degree $\lqe n$.  Then $U(\g_{-})=\cup U(\g_{-})_{n}$.  We give $U(\g_{-})_{n}$ the discrete topology and then $U(\g_{-})=M_{+}$ with the inductive limit topology, i.e. in this case the discrete topology.  
Let $M_{+}^{*}$ be the superspace of all continuous linear functionals on $M_{+}$.   For any $n \geq 0$ let $U(\g_{-})_{n}$ be the elements of $U(\g_{-})$ with degree $\leq n$.  Then $M_{+}^{*}$ is the projective limit of $U(\g_{-})_{n}^{*}$.  By giving $U(\g_{-})_{n}^{*}$ the discrete topology, the superspace $M_{+}^{*}$ inherits a natural structure of a topological superspace.

Let $i_{+}^{*}:M_{+}^{*}\ttp M_{+}^{*}\rightarrow M_{+}^{*}$ be the map defined by
$$i_{+}^{*}(f\otimes g)(x) := (f\otimes g)i_{+}(x)$$
for $f,g \in M_{+}^{*}$ and $x\in M_{+}$.  By definition of the topology on $M_{+}^{*}$ the map $i_{+}^{*}$ is continuous.  Therefore, $i_{+}^{*}$ extends to a morphism $i_{+}^{*}:M_{+}^{*} \ttp M_{+}^{*} \rightarrow M_{+}^{*}$.  The proof of Lemma 8.3 \cite{EK1} implies that the m $i^{*}_{+}$ associative, i.e. $ i^{*}_{+} \circ (i^{*}_{+}\otimes 1) \Phi^{-1}=  i^{*}_{+} \circ(1 \otimes i^{*}_{+}) $.

\subsection{The Tensor Functor $\FF$}

Define the functor $\FF:\Mt \rightarrow \At$ as
\begin{equation}
\label{E:DefF}
\FF(V)=\Hom_{\M}(M_{-},M^{*}_{+} \ttp V).
\end{equation}
From the following Lemma \ref{L:FV} we have that $\FF:V \rightarrow V[[h]]$ where $V[[h]]$ is the topologically free $k[[h]]$-module associated to the graded vector space underlying $V$.
For any $V\in \M$ define 
$$\Psi_{V}: \Hom_{\g}(M_{-},M^{*}_{+}\ttp V) \rightarrow V$$ by $f \rightarrow (1_{+}\otimes 1)f(1_{-})$ where $(1_{+}\otimes 1)(g\otimes v):=g(1_{+})v$ for $g \in M^{*}_{+}$ and $v\in V$. 

\begin{lem}\label{L:FV}
 $\Psi_{V}$ is a even vector superspace isomorphism.
\end{lem}

\begin{pf}
The proof of the lemma follows from checking that the isomorphisms of \cite[Lemma 8.1]{EK1} preserve the $\Z_{2}$-grading.  We define the these isomorphisms and see that they are even.
 
 By Frobenius reciprocity the following maps 
$$\Hom_{\g}(M_{-},M_{+}^{*}\ttp V)\rightarrow  (M^{*}_{+}\ttp V)^{\g_{-}}  \text{ given by } f \mapsto f(1_{-}), $$ 
$$ \Hom_{\g_{-}}(M_{+},V)\rightarrow V \text{  given by }  f \mapsto f(1_{+})$$
are isomorphism of topological vector spaces.  
Let $$ (M^{*}_{+}\ttp V)^{\g_{-}} \rightarrow  \Hom_{\g_{-}}(M_{+},V)$$ be the map given by $f\otimes x \mapsto f_{x}$, where $f_{x}(y):= (-1)^{\p y \p x}f(y)x$.  This map is an isomorphism of topological vector spaces (see \cite[Proof of Lemma 8.1]{EK1}) where $ \Hom_{\g_{-}}(M_{+},V)$ has the weak topology.   Also by definition all of these maps are even homomorphisms of superalgebras.  Composing the above maps we have the desired isomorphism 
 $$\Hom_{\g}(M_{-},M_{+}^{*}\ttp V)\rightarrow   V  \text{ which is given by } f \mapsto (1_{+}\otimes 1)f(1_{-}). $$ 
\end{pf}

\begin{defn}
For each pair $V,W \in \M$ define $\FJ_{V,W}: \FF(V)\atp \FF(W) \rightarrow \FF(V\ttp W)$ by 
$$\FJ_{V,W}(v\otimes w)=(i^{*}_{+}\otimes 1 \otimes 1) \circ \Phi^{-1}_{1,2,34} \circ (1\otimes \Phi_{2,3,4}) \circ  \beta^{-1}_{23} \circ (1\otimes \Phi^{-1}_{2,3,4}) \circ \Phi_{1,2,34} \circ (v \otimes w)\circ i_{-}.$$
\end{defn}
%When it is clear we will not write the associative isomorphism.  Using this convention:
%$$\FJ_{V,W}(v\otimes w)=(i^{*}_{+}\otimes 1 \otimes 1) \circ  \beta^{-1}_{23} \circ(v \otimes w)\circ i_{-}.$$
\begin{thm}\label{T:J}
The collection $(\FJ_{V,W})_{V,W \in \M}$ defines a tensor structure on $\FF$, i.e. $\FF$ is a tensor functor.
\end{thm}
\begin{pf}
Using the facts that $i_{-}$ is coassociative and $i^{*}_{+}$ is associative the proof follows exactly in the same way as the universal or pictorial proof of Proposition 19.1 \cite{ES}.
\end{pf}

Let $\End(\FF)$ be the endomorphisms of $\FF$ (see \S \ref{SS:EndofF}).  Using $\Psi_{V}$ to identify $\FF(V)$ and $V[[h]]$ the proof of Lemma \ref{L:UHF} shows that there exists a canonical even superalgebra isomorphism 
\begin{equation}
\label{E:U[[h]]isoEndF2}
\theta : \U[[h]] \rightarrow  \End(\FF) \text{ given by } x \mapsto x_{|V}
\end{equation}
where $x_{|V}$ is $x$ acting on the $\U[[h]]$-module $V[[h]]$.  We use this isomorphism is to identify $\End(\FF)$ and $\U[[h]]$.

Next we will define an element  $\FJ \in \U^{\ttp 2}[[h]]$ whose action on $V[[h]] \atp W[[h]]$ determines $\FJ_{V,W}$.  Recall the isomorphism $\Psi_{V}: \Hom_{\g}(M_{-},M^{*}_{+}\ttp V) \rightarrow V$ of Lemma \ref{L:FV}.  Let $\phi : M_{-}\rightarrow M^{*}_{+}\ttp \U$ be the even morphism given by $\phi=\Psi^{-1}_{\U}(1)$.  Given $g \in \Hom_{\g}(V,W)$ denote the map $\hat{g}:V[[h]] \rightarrow W[[h]]$ by $\sum v_{i}h^{i} \rightarrow \sum g(v_{i})h^{i}$.  

Define the element $\FJ \in \U^{\ttp 2}[[h]]$ by
$$\FJ=(1_{+}\otimes 1)\left(( i^{*}_{+}\otimes 1 \otimes 1)  \Phi^{-1}_{1,2,34}  (1\otimes \Phi_{2,3,4})    \beta^{-1}_{23}  (1\otimes \Phi^{-1}_{2,3,4})  \Phi_{1,2,34}  \left(y \right)\right)$$
where $y = \phi(1_{-}) \otimes  \phi(1_{-} )$.

The following lemma shows that the map $\FJ_{V,W}$ is determined by the element $\FJ$.
\begin{lem}
\label{L:ActionOfJ}
 Let $\theta$ be the isomorphism given in (\ref{E:U[[h]]isoEndF2}).  Then $\theta(\OJ)=\OJ_{V,W}$, i.e. 
\begin{eqnarray}
\label{E:AJ}
\FJ(v\otimes w) &= &\hat{\Psi}_{V\otimes W}\big( \FJ_{V,W}(\hat{\Psi}^{-1}_{V}v \otimes \hat{\Psi}^{-1}_{W}w)\big)
\end{eqnarray}
for all $v \in V[[h]]$ and $w \in W[[h]]$.
\end{lem}
%
% Proof of Lemma about Action of J 
%
\begin{pf} 
By (\ref{E:Uhmodules}) it is enough to check that (\ref{E:AJ}) holds for $v\in V$ and $w\in W$. 
Let $\phi(1_{-})$ be the tensor $\sum f_{i}\otimes x_{i} \in (M^{*}_{+}\ttp \U)^{\g_{-}}$, then we have $\sum f_{i}(1_{+})x_{i}=1$.
Let $v \in V$ and consider the following calculation:
$$(1_{+}\otimes 1)(\sum f_{i} \otimes (x_{i}v))=\sum f_{i}(1_{+})x_{i}v=v. $$
The above shows that $\Psi_{V}^{-1}(v)(1_{-})=\sum f_{i} \otimes (x_{i}v)$.  

Let $\Fvtheta \in \U^{\ttp 4}[[h]]$ be given by:
\begin{eqnarray*}
\Ftheta & := & \Phi^{-1}_{1,2,34}  (1\otimes \Phi_{2,3,4})  \beta^{-1}_{23} (1\otimes \Phi^{-1}_{2,3,4})  \Phi_{1,2,34}  
% & = &  \Phi^{-1}_{1,2,34}  (1\otimes \Phi_{2,3,4})    \flip_{23} e^{\frac{h t_{23}}{2}}  (1\otimes \Phi^{-1}_{2,3,4})  \Phi_{1,2,34} \\
 %&=&  \Phi^{-1}_{1,2,34}  (1\otimes \Phi_{2,3,4})  e^{\frac{h t_{23}}{2}}  (1\otimes \Phi^{-1}_{3,2,4}) \Phi_{1,3,24}. 
\end{eqnarray*}
We use $\Fvtheta$ to simplify notation.
Represent $\Fvtheta=\sum_{i} \Ftheta_{i}h^{i}$ where $\Ftheta_{i}= \sum_{k}  \Ftheta_{i}^{ 1k}\otimes \Ftheta_{i}^{2 k}\otimes \Ftheta_{i}^{3 k}\otimes \Ftheta_{i}^{ 4k}$.
Evaluating the right side of (\ref{E:AJ}) with $v  \in V$ and $w \in W$, we have:
\begin{align*}
\hat{\Psi}_{V\ttp W}\big( \FJ_{V,W} & (\hat{\Psi}^{-1}_{V}v \otimes \hat{\Psi}^{-1}_{W}w)\big)  = \\ 
& = (1_{+}\otimes 1)\big[i_{+}^{*}\circ \Fvtheta \circ \flip \circ \sum_{j} f_{j} \otimes (x_{j}v) \otimes \sum_{l} f_{l} \otimes (x_{l}w) \big]\\
 & =  (1_{+}\otimes 1)\bigg[ i_{+}^{*}\circ \sum_{i}\Ftheta_{i}h^{i} \bigg( \sum_{j,l} (-1)^{\overline{x_{j}v}\p{f_{l}}} f_{j} \otimes f_{l} \otimes (x_{j}v) \otimes  (x_{l}w) \bigg) \bigg] \\
 & = \bigg[ \sum_{i,j,l,k}h^{i}(-1)^{A}f_{j}(\Ftheta_{i}^{1k})f_{l}(\Ftheta_{i}^{ 2k})\Ftheta_{i}^{3k}x_{j} \otimes  \Ftheta_{i}^{4 k}x_{l}\bigg]v\otimes w
\end{align*}
where $A=\overline{x_{j}v}\p{f_{l}}+ \p{f_{j}}(  \overline{ \Ftheta_{i}^{2 k} } + \overline{ \Ftheta_{i}^{ 3k} }   + \overline{ \Ftheta_{i}^{4 k} } )    + \p{f_{l}} (  \overline{ \Ftheta_{i}^{3 k} }+  \overline{ \Ftheta_{i}^{4 k} })  +  \overline{x_{j}v} \: \overline{ \Ftheta_{i}^{4 k} }  +  \p{f_{j}}  \overline{ \Ftheta_{i}^{ 1k} }  + \p{f_{l}}    \overline{ \Ftheta_{i}^{ 2k} }    + ( \overline{ \Ftheta_{i}^{ 4k} x_{l}}) \p v $.

Similarly evaluating the left side of (\ref{E:AJ}) we have:
 \begin{eqnarray*}
\FJ(v\otimes w) & = & \big[ (1_{+} \otimes 1)(i^{*}_{+}\circ \Ftheta \circ \flip_{23 } \circ \phi \otimes \phi (1_{-}\otimes 1_{-})) \big] v\otimes w\\
 & = & \bigg[ (1_{+}\otimes 1)\big(i_{+}^{*}\circ \sum_{i}\Fvtheta_{i}h^{i} ( \sum_{j,l} (-1)^{\overline{x_{j}}\p{f_{l}}} f_{j} \otimes f_{l} \otimes x_{j} \otimes  x_{l} )\big)\bigg]v\otimes w\\
 & = & \bigg[\sum_{i,j,k,l}(-1)^{B}h^{i}f_{j}(\Ftheta_{i}^{1k})f_{l}(\Ftheta_{i}^{ 2k})\Ftheta_{i}^{3k}x_{j} \otimes  \Ftheta_{i}^{4 k}x_{l}\bigg]v\otimes w
\end{eqnarray*}
where $B=A-\p{v}\p{f_{l}}- \p{v} \: \overline{ \Ftheta_{i}^{4 k} }  -( \overline{ \Ftheta_{i}^{4 k}. x_{l}}) \p v =A-\p{v}\p{f_{l}}- \p{v} \: \overline{ \Ftheta_{i}^{4 k} }  -\overline{\Ftheta_{i}^{4 k}} \p v - \p{x_{l}}\p v= A$.  The last equality follows from the fact that $\sum f_{l}\otimes x_{l} $ is even, i.e. $\p{f_{l}}=\p{x_{l}}$. Thus we have showed that (\ref{E:AJ}) holds, completing the proof.
\end{pf}

%\subsection{The quantized enveloping algebra $U_{h}(\g)$} 
\subsection{The quantization of the double $\g= D(\g_{+})$.}

As in \S \ref{S:NQof}, we will now define a Hopf superalgebra structure on $\U [[h]]$ and show it is a quantization of $\g$.  After replacing $\OJ$ with $\FJ$ equations (\ref{E:DefCoproduct}) and (\ref{E:DefAntip}) define a Hopf superalgebra structure on $\U[[h]]$.   Let $\FH$ be this Hopf superalgebra.  

% Theorem  \\FH is quant of \g	

\begin{thm}
\label{T:Qofg2}
	$\FH$ is a quantization of the Lie superbialgebra \g.
\end{thm}

\begin{pf}
We need the following lemma.
\begin{lem} 
\label{L:Jmod2}

$\FJ \equiv 1 + \frac{rh}{2} \mod h^{2}$

\end{lem}

\begin{pf}[Proof of Lemma \ref{L:Jmod2}]

Since $\Phi \equiv 1 \mod h^{2}$, we have
\begin{align}
   \FJ & \equiv (1_{+}\otimes 1)[i^{*}_{+}(1 - \tfrac{t_{23}h}{2})\flip_{23}(\phi(1_{-}) \otimes \phi(1_{-}))]   \mod h^{2}  \notag  \\
    &   \equiv (1_{+}\otimes 1)\bigg[ i^{*}_{+}\left(1 - \tfrac{(r_{23}+\flip_{23}r_{23})h}{2}\right)\bigg(\sum_{i,j} (-1)^{\p{f_{j}}\p{x}_{i}}f_{i}\otimes f_{j} \otimes x_{i} \otimes x_{j} \bigg)\bigg] \mod h^{2} \notag \\
    &  \equiv \sum_{i,j} (-1)^{\p{f_{j}}\p{x}_{i}} f_{i}(1_{+})f_{j}(1_{+})x_{i}\otimes x_{j}  \notag   \\
    & \hspace{25pt}  -  \tfrac h2 \sum_{i,j,k} (-1)^{  \p{f_{j}}\p{x}_{i}  +     \p{f_{j}}\p{m}_{k}  + 1  +  \p{f_{j}}\p{p}_{k}   +  \p{p}_{k} \, \p{m}_{k}  }  f_{i}(1_{+})f_{j}(m_{k}1_{+})p_{k}x_{i}\otimes x_{j}   \mod h^{2} \notag \\
    & \equiv 1 -   \tfrac h2 \sum_{j,k} (-1)^{     \p{f_{j}}\p{m}_{k}  + 1  +  \p{f_{j}}\p{p}_{k}   +     \p{p}_{k} \, \p{m}_{k}   }  f_{j}(m_{k})p_{k}\otimes x_{j}   \mod h^{2}  \label{E:Jmod4}  \\
    &  \equiv 1 +  \tfrac h2 \sum_{j,k} p_{k}\otimes (-1)^{ \p{m}_{k}} f_{j}(m_{k})x_{j}   \mod h^{2}  \label{E:Jmod5}  \\
    &  \equiv 1 + \tfrac{hr}{2} \mod h^{2}  \label{E:Jmod6} 
\end{align}
where $r=\sum_{k } p_{k}\otimes m_{k} \in \g_{+}\otimes \g_{-}$ is the canonical element of $D(\g_{+})$ defined in \S\ref{S:Lbialg}.  The first three equivalences follow by definition.  Equivalence (\ref{E:Jmod4}) follows from the facts: $\p{x}_{i}=\p{f_{i}}$; if $f_{i}$ is odd, then $f_{i}(1_{+})=0$ and $\sum f_{i}(1_{+})x_{i}=1$.  Equivalence (\ref{E:Jmod5})  follows from the fact that $r$ is even,  (\ref{E:Jmod6}) hold because of the identity $\sum (-1)^{\p{m}_{k}} f_{i}(m_{k})x_{i}=m_{k}$ (which follows from $\p{m}_{k} \neq \p x_{i}$ implies $f_{i}(m_{k})=0$).    Thus we have proven the Lemma.
\end{pf}
From Lemma \ref{L:Jmod2} it follows that the equivalences (\ref{E:Qofg1}) and (\ref{E:Qofg2}) hold for all $x \in \g \subset \FH$.  The theorem is proved.
\end{pf}

%   The quantization of the Lie bialgebra $\g_{+}$.

\subsection{The quantization of the Lie bialgebra $\g_{+}$.}

	As in \S\ref{S:NQof} we use $\theta$ to identify $\FH$ and $\End(\FF)$.  In this subsection we will construct a Hopf subalgebra of $\FH$ denoted $\FUh(\g_{+})$ which will be a quantization of the Lie bialgebra $\g_{+}$.

	Consider the even superspace homomorphism $i:\FF(M_{-})\rightarrow  \FH$ given by $x\mapsto i(x)$ where 
$$i(x)v=(-1)^{\p x \p v}(i^{*}_{+}  \otimes 1)\circ (1 \otimes v) \circ x$$
for any $V \in \M$ and $x \in \FF(M_{-}) $ and $ v \in \FF(V)$.  Next we will show that the map $i$ is injective.  Recall the isomorphism $\hat \Psi : \FF(M_{-})  \rightarrow U(\g_{+})[[h]]$ given by $f \mapsto (1_{+}\otimes 1) f(1_{-})$.  For $x \in U(\g_{+})[[h]]$ let $f_{x}:=\hat \Psi^{-1}(x)$.  Now for $x \in U(\g_{+})$ and $v\in \FF(V)$ we have $\hat \Psi (i(f_{x} )v) \equiv x\hat \Psi(v) \mod h$.  Therefore, $i$ is injective.  
Set $\FUhp =i(\FF(M_{-}))$.   
 
%     theorem}	\FUhp is a quantization of the Lie superbialgebra $\g_{+}$
%
\begin{thm}
		$\FUhp$ is a quantization of the Lie superbialgebra $\g_{+}$.
\end{thm}
\begin{pf}

The following lemma implies that $\FUhp$ is a sub-superbialgebra of $\FH$. 

%Lemma    \FUhp is closed under multiplication and coproduct in $\Uh$
%
\begin{lem}
\label{L:UhpClosed}
	$\FUhp$ is closed under multiplication and coproduct in $\FH$.
\end{lem}
%
%
%Proof of Lemma    \FUhp is closed under multiplication and coproduct in $\Uh$
%
%
\begin{pf}
	For $x,y \in \FF(M_{-})$ and $v \in \FF(V)$ the associativity of $i_{+}^{*}$ and relation (\ref{R:Pent}) imply (to simplify notation set $w=(1 \otimes y) x$)
\begin{align*}
   i(x) \circ i(y)v & =(-1)^{\p x \p y +  \p v \p x  + \p y \p v}(i^{*}_{+}  \otimes  1)\Phi^{-1} (1\otimes i^{*}_{+}  \otimes  1)\Phi^{-1}_{2,3,4} (1 \otimes 1 \otimes v)(1 \otimes y) x  \\
   & =(-1)^{\p x \p y +  \p v \p x  + \p y \p v}(i^{*}_{+}  \otimes 1) ( i^{*}_{+}  \otimes 1  \otimes  1) \Phi^{-1}_{1,2,3} \Phi^{-1}_{1,23,4}  \Phi^{-1}_{2,3,4}  (1 \otimes 1 \otimes v) w  \\
      & =(-1)^{\p x \p y +  \p v \p x  + \p y \p v}(i^{*}_{+}   \otimes 1) ( i^{*}_{+}  \otimes 1 \otimes  1) \Phi_{12,3,4}^{-1}  \Phi_{1,2,34}^{-1}  (1 \otimes 1 \otimes v) w \\
%          & =(-1)^{\p x \p y +  \p v \p x  + \p y \p v}(i^{*}_{+}  \otimes 1) \Phi^{-1} ( i^{*}_{+}  \otimes 1 \otimes 1)  (1 \otimes 1 \otimes v) \Phi^{-1} (1 \otimes y) x  \\
         & =(-1)^{\p x \p y +  \p v \p x  + \p y \p v}(i^{*}_{+}  \otimes  1)\Phi^{-1}  (1 \otimes v) (i^{*}_{+}  \otimes  1) \Phi^{-1}  (1 \otimes y)  x  \\
    &  =(-1)^{ \p v \p x  + \p y \p v}(i^{*}_{+}  \otimes  1)\Phi^{-1}  (1 \otimes v) z \\
    & =  i(z)v
\end{align*}
where $z= (-1)^{\p x \p y}(i^{*}_{+}  \otimes  1)\circ \Phi^{-1} \circ (1 \otimes y) \circ x $.

	Following the proof  in \cite{EK1} Chapter 9 we have
$$ \Delta (i(x)) = (i \otimes i)(\FJ^{-1}_{M_{-},M_{-}}(1 \otimes i_{-}) \circ x)$$
which completes the proof of the Lemma.
\end{pf}
Next we show that $\FUhp$ is a Hopf superalgebra.  Consider the even superspace isomorphism 
$$\mu :U(\g_{+})[[h]]\rightarrow \FUhp \text{ given by } x \mapsto i(f_{x})$$
where $f_{x}:=\hat \Psi^{-1}(x)$.  For $x,y \in U(\g_{+})[[h]]$ we have
\begin{align*}
   i(f_{x})\circ i(f_{y}) & \equiv  i \left( (-1)^{\p x \p y} (i^{*}_{+} \otimes  1) \circ (1 \otimes f_{y}) \circ f_{x} \right) & \mod h^{2} \\
    &  \equiv i(f_{xy})   &   \mod h^{2}
\end{align*}  
i.e. $\mu (x) \circ \mu (y) \equiv \mu (xy)\mod h^{2}$.
Similarly, we have $(\mu \otimes \mu) \D{x} \equiv \D{\mu (x)} \mod h $.  Therefore, $\FUhp / h\FUhp $ is isomorphic to $U(\g_{+})$ as a superbialgebra.  This implies that $\FUhp$ has a Hopf superalgebra structure. 

To finish the proof we need to show that the equivalence (\ref{R:ConditionDelta}) holds.
%	Next we will show that the inclusion $U(\g_{+}) \rightarrow \Uh$ is nice module $h^{2}$.  
Recall the isomorphism $\theta : \U[[h]] \rightarrow  \FH$ given in (\ref{E:U[[h]]isoEndF2}).   Then we have
 \begin{equation}
\label{E:muEquivTheta}
  \mu (x) \equiv \theta (x) \mod h^{2}
\end{equation}
for all $x \in U(\g_{+})$.  In other words, the image of $U(\g_{+}) $ in  $\FUhp$ and $\FH$ is equal modulo $h^{2}$.  From Theorem \ref{T:Qofg2} we have that the equivalences (\ref{E:Qofg1}) and (\ref{E:Qofg2}) hold for all $x \in \g \subset \FH$.  Combining the last statement with (\ref{E:muEquivTheta}) and the fact that $\g_{+}$ is a Lie sub-superbialgebra of $\g$ we have that the equivalences (\ref{E:Qofg1}) and (\ref{E:Qofg2}) hold for all $x \in \g_{+} \subset \FUh(\g_{+})$. Thus, $\FUhp$ is a quantization of $\g_{+}$. 
\end{pf}

\begin{thm}\label{T:twoQareIso}
 Let $\g_{+}$ be a finite dimensional Lie superbialgebra.  The quantization of $\g_{+}$ constructed in \S \ref{S:NQof} isomorphic to the be the quantization of $\g_{+}$ constructed in this section, i.e. $\OUh(\g_{+}) \cong \FUh(\g_{+})$.
\end{thm}
\begin{pf} Let $\g$ be the double of $\g_{+}$ and let $\tMt$ be the category of discrete $\g$-modules.  Consider the functor $\tF:\Mt\rightarrow \At$ given by
$$\tF(V)=\End_{\Mt}(M_{+}\ttp M_{-},V).$$
 By definition $\End(\tF |_{\tMt})$ is the quantization $\OH$ of the double $\g$, defined in \S \ref{S:NQof}.  Since $\End(\tF) $ and $\OH$ are both isomorphic to $\U[[h]]$, we have that the morphism $\zeta :\End(\tF) \rightarrow \OH$ given by the restriction of $\Mt$ to $\tMt$, is an isomorphism of Hopf superalgebras.   

Let $\chi: \tF \rightarrow \FF$ be the natural transformation of functors given by $\chi_{V}(v)=(1\otimes v)\circ (\sigma \otimes 1)$ where $\sigma$ is the canonical element in $\Hom_{\Mt}(k,M_{+}^{* }\otimes M_{+})$.  Using the properties of the braiding $\beta$ one can follow the proof of Proposition 9.7 in \cite{EK1} to show that $\chi$ is a natural isomorphism of tensor functors.  Therefore, $\chi$ induces an isomorphism between the Hopf superalgebras $\End(\tF)$ and $\End(\FF)$.  Composing this isomorphism with $\zeta^{-1}$ we have an isomorphism of Hopf superalgebras $\kappa: \OH \rightarrow \FH$.   By construction the image of the  restriction of $\kappa$ to the Hopf sub-superalgebra $\OUhp$ is $\FUh(\g_{+})$.  In other words, $\kappa|_{\OUhp}:\OUhp \rightarrow \FUhp$ is an isomorphism of Hopf superalgebras.   
\end{pf}
Using the isomorphism $\kappa$ ($\kappa|_{\OUhp}$) given in the proof of theorem \ref{T:twoQareIso} we will identify $\OH$ and $\FH$ (resp. $\OUhp$ and $\FUhp$).  From this point on, we will make no distinctions  
between $\OH$ and $\FH$ or $\OUhp$ and $\FUhp$.   We call $\OUhp$ the Etingof-Kazhdan quantization of $\g_{+}$.

%  Section     Functoriality of quantizations    

% marker
\section{Functoriality of the quantizations}\label{S:FunctQ}
In this section we show that the quantizations \ref{S:2QUofg} and \ref{S:QUofquasi} are functorial. Then we use this to show that the quantization commutes with taking the double.

Let $LSBA(k)$ be the category of finite dimensional Lie superbialgebra over $k$ and let $QUES(K)$ be the category of QUE superalgebra over $K=k[[h]]$.

%  Theorem  functor LSBA to QUES

\begin{thm}
There exists a functor from $LSBA(k)$ to $QUES(K)$ such that $\liea \in LSBA(k)$ is mapped to $U_{h}(\liea)$ which is the quantization defined in \S \ref{S:2QUofg}. 
\end{thm}

\begin{pf}
Once one accounts for the necessary signs, the proof is identical to the classical case (c.f. Theorem 10.1 and 10.2, \cite{EK1}).
\end{pf}

%Note the proof also follows from the work of Enriquez \cite{En}, where he quantizes Lie bialgebras at the universal level using cohomology computations.  Moreover, Enriquez results show that his twist is gauge equivalent to our twist $J$, and therefore our constructions are equivalent.  In particular our quantization is functorial.     

%   Def of   QTLSBA(k)$ be the category

Let $QTLSBA(k)$ be the category of quasitriangular Lie superbialgebra over $k$ and let $QTQUES(K)$ be the category of quasitriangular QUE superalgebra over $K=k[[h]]$.

%  Theorem  functor LSBA to QUES

\begin{thm}
There exists a functor from $QTLSBA(k)$ to $QTQUES(K)$ such that $(\g,r) \in QTLSBA(k)$ is mapped to $(U_{h}^{qt}(\g),R)$ which is the quantization defined in \S \ref{S:QUofquasi}. 
\end{thm}

\begin{pf}
 The proof is a consequence of Theorem 1.2 (ii) of \cite{EK2}.  Theorem 1.2 (ii) states that there is a ``universal quantization functor''  from the cyclic category of quasitriangular Hopf algebras to the closure cyclic category of quasitriangular Lie bialgebras (see \cite{EK2}).  By considering linear algebraic structures in the symmetric tensor category of superspaces this ``universal quantization functor'' gives rises to functor from $QTLSBA(k)$ to $QTQUES(K)$ with the desired properties.   
%     We will show that the proof is a consequence of Theorem 1.2 (ii) in \cite{EK2}.   Let $\veccat$ ($\supcat$) be the symmetric tensor category of vector spaces (resp. superspaces).  As in \cite{EK2}, let $\qtlba$ ($\qtha$) be the cyclic category of quasitriangular Lie bialgebras (quasitriangular Hopf algebras).  Loosely speaking, $\qtlba$ is a tensor category generated by one element (the canonical quasitriangular Lie bialgebras) such that the morphism are generated by morphisms reflecting the structure of a quasitriangular Lie bialgebras.  Theorem 1.2 (ii) states that there exists a ``universal quantization functor'' 
%$$Q :  \qtha\rightarrow \qtlba$$
%
%
%
%Follows from Theorem 1.2 (ii) in \cite{EK2}, which states that there is a ``universal quantization functor'' universal construction from the cyclic category generated by the canonical qua tri H algebras to the clo of cyclic category generated by the canonical q t Lie bialgebra such that the functor induced from     this theorem holds for the cyclic for any symmetric tensor category.   In particular, the category of superspaces where flip..
\end{pf}
%
%
%
%Note the proof of Theorem 1.2 (i) in \cite{EK2}, requires elements check that it is a sub-algebra.  May be we should say that is why we did not use this theorem to define \Hgff  
%
%
%
%

Next we use the functoriality to prove the following theorem which first appeared in \cite{EK1} for the non-super case.

\begin{thm}\label{T:QtQisIsm}
Let $\g_{+}$ be a finite dimensional quasitriangular Lie superbialgebra.  Then the quantization of the quasitriangular Lie superbialgebra $\g_{+}$ constructed in \S \ref{S:QUofquasi} is isomorphic to the quantization of the Lie superbialgebra $\g_{+}$ of \S \ref{S:2QUofg}, i.e. $U_{h}^{qt}(\g_{+}) \cong U_{h}(\g_{+})$ as Hopf algebras.
\end{thm}

\begin{pf}
To prove the theorem we need the following lemma (which first appeared in \cite{RSTS} for the non-super case).

% Lemma  morphism 

\begin{lem} 
Let $(\g_{+},r)$ be quasitriangular Lie superbialgebra and $\g=D(\g_{+})$ be its double.  Then there exist a quasitriangular Lie superbialgebra morphism $ \g \rightarrow \g_{+}$, which is the identity when restricted to $\g_{+}$.
\end{lem}
\begin{pf} 
Let $\upsilon : \g=\g_{+}\oplus \g_{+}^{*} \rightarrow \g_{+}$ be the linear map given by 
$$\upsilon(x+f)=-x-(1\otimes f)r$$
for $x \in \g_{+}$ and $f \in \g_{+}^{*}$.  We will show that $\upsilon $ is a quasitriangular Lie superbialgebra morphism.  

First we show it is a Lie superbialgebra morphism, i.e. $\upsilon([a,b])=[\upsilon(a),\upsilon(b)]$ for all $a,b \in \g$.  This is clear if $a,b \in \g_{+}$.  Recall the definition of the bracket on the double given in (\ref{R:BracketDouble}).  Then if $x \in \g_{+}$ and $ f \in \g_{+}^{*}$ we have 
\begin{align*}
\label{}
   \upsilon([x,f]) &= \upsilon \left((ad^{*}x)f-(-1)^{\p x \p f}(1 \otimes f)\delta(x) \right)  \\
    &  = -(-1)^{\p x \p f}(1 \otimes f)[1 \otimes x, r] + (-1)^{\p x \p f}(1 \otimes f)[x \otimes 1 + 1 \otimes x,r] \\
    & =  (-1)^{\p x \p f}(1 \otimes f)[x \otimes 1, r] \\
    & =  [x, (1 \otimes f)r]  \\
    & = [ \upsilon (x), \upsilon(f)]
\end{align*}
Note that $(ad^{*}x)f  $ is the linear functional $y \mapsto (-1)^{\p x \p f}f \circ [x,y]$.  Similarly, one shows that $\upsilon([f,g])=[\upsilon(f),\upsilon(g)] $ for $f,g \in \g_{+}^{*}$.

Finally, we need to show that $\upsilon $ is a quasitriangular Lie superbialgebra morphism, i.e. preserves the r-matrix.  Let $\check{r}$ be the r-matrix of $\g$.  Choose a basis $x_{i}$ for $\g_{+}$ and let $f_{i}$ be the dual basis of $\g_{+}^{*}$, then $\check{r}= \sum x_{i}\otimes f_{i}$.  Therefore we have
$$(\upsilon \otimes \upsilon)(\check{r})=\sum \upsilon(x_{i}) \otimes \upsilon (f_{i})= \sum x_{i}\otimes(1\otimes f_{i})r =r.$$ 
Thus $\upsilon$ is the desired morphism.
\end{pf}

Now we prove the theorem.  Recall that by construction $U_{h}(\g_{+})$ is a subalgebra of $\OH$.  From the Lemma we have $\upsilon : \g \rightarrow \g_{+}  $ such that $\upsilon |_{\g_{+}}=id_{\g_{+}}$.  The functoriality of the quantization implies that $\upsilon $ induces a morphism of QTQUE superalgebras $\Uqt \rightarrow U_{h}^{qt}(\g_{+})$.  Restricting this morphism to the subalgebra $U_{h}(\g_{+})$ we have a morphism $U_{h}(\g_{+}) \rightarrow U_{h}^{qt}(\g_{+})$, which is a isomorphism since it is modulo $h$.
\end{pf}

We end this section with the following theorem.  

%  Thm        The quant commutes with taking the double.

\begin{thm}
The quantization of a finite dimensional Lie superbialgebra $\g_{+}$ commutes with taking the double, i.e. $D(\Uhp)\cong U_{h}(D(\g_{+}))$ (for the definitions of the doubles see \S \ref{S:Lbialg} and Proposition \ref{P:DoubleQUE}). 
\end{thm}
\begin{pf} 
From Corollary \ref{C:qd} we have $D(\Uhp)\cong \OH$, where $\OH$ is the quantization of $D(\g_{+})$ constructed in \S \ref{S:NQof}.  By construction $U_{h}^{qt}(D(\g_{+}))=\OH$, where $U_{h}^{qt}(D(\g_{+}))$ is quantization of $D(\g_{+})$ given is \S \ref{S:QUofquasi}.  By Theorem \ref{T:QtQisIsm} we have $ U_{h}^{qt}(D(\g_{+}))\cong U_{h}(D(\g_{+}))$.  Combining the above isomorphism we have the desired result. 
\end{pf}

%           SECTION      E-K Quantization of the Lie superbialgebra type A-G

%marker
\section{The Etingof-Kazhdan quantization of Lie superalgebras of type A-G}\label{S:QofglConclution}
In this section we will show that, for Lie superalgebras of type A-G, the E-K quantization is isomorphic to the Drinfeld-Jimbo quantization.  We follow \cite{EK6} which proves the result for generalized Kac-Moody algebras.  However, we must take the new quantum Serre-type relations into consideration.  As in \cite{EK6} we will show that the E-K quantization is given by the desired generators and relations.  In particular, we extend results of Lusztig \cite{L} to the setting Lie superalgebras of type A-G and check directly that the new quantum Serre-type relations are in the kernel of the appropriate bilinear form.    

Here we recall some notation from \S \ref{SS:g} and \S\ref{SS:glbialg}.  Let $\g$ be a Lie superalgebra of type A-G.  Let $\Phi=\{\alpha_{1},...,\alpha_{s}\}$ be a simple root system with at most one odd root and let $(A,\tau)$ be the corresponding Cartan matrix where $\tau=\{m\}$  or $\tau=\emptyset$.  Let $d_{1},\dots,d_{s}$ be the nonzero numbers such that $d_{i}a_{ij}=d_{j}a_{ji}$ and $d_{1}=1$.  Let $(,)$ be the unique non-degenerate supersymmetric invariant bilinear form on $\g$.  By rescaling if necessary we may assume that the restriction of $(,)$ to $\h$ is determined by $(a,h_{i})=d_{i}^{-1}\alpha_{i}(a)$ for all $a\in \h$ and $i\in I=\{1,...,s\}$.

Let $\ghat$ be the Lie superalgebra generated by $e_{i}, f_{i}$ and $h_{i}$ for $i\in I$ satisfying (\ref{R:LieSuperalg}) where all generators are even expect for $e_{t}$ and $f_{t}$ when $t\in \tau$ which are odd.  Let $\borelhat_{\pm}$ be the Borel sub-superalgebra of $\ghat$ generated by $e_{i}, h_{i}$ and $f_{i}, h_{i}$, respectively.  Let $q=  h/2$.

%    Subsection     The quantized enveloping superalgebra $\Uhbhat$

\subsection{Generators and relations for $\Uhbhat$}

\begin{thm}\label{T:GRofUhbhat}
The quantized universal enveloping superalgebra $\Uhbhat$ is isomorphic to the quantized enveloping superalgebra $\Uhatp$ generated over $\C[[h]]$ by the elements $e_{i}, h_{i},\: i\in I$ (where all generators are even expect for $e_{t}$, $t\in \tau$ which is odd) satisfying the relations
\begin{align*}
\label{}
    [h_{i},h_{j}]=&0, & [h_{i},e_{j}]=&a_{ij}e_{j},
\end{align*}  
with coproduct
\begin{align*}
\label{}
    \D{h_{i}}=&1\otimes h_{i} + h_{i} \otimes 1, & \D{e_{i}}=&e_{i}\otimes q^{d_{i}h_{i}}+ 1\otimes e_{i},
\end{align*}   
for all $i,j\in I$.
\end{thm}

The theorem follows from the following two lemmas.  

%   LEMMA  UQES alg Uhbhat iso to ...

\begin{lem}\label{L:UQESbhat}
The universal quantized enveloping superalgebra $\Uhbhat$ is isomorphic to the quantized enveloping superalgebra generated over $\C[[h]]$ by the elements $e_{i}, h_{i}, \: i\in I$ (where all generators are even expect for $e_{t}$, $t\in \tau$ which is odd) satisfying the relations
\begin{align*}
\label{}
      [h_{i},h_{j}]=&0, & [h_{i},e_{j}]=&a_{ij}e_{j},
\end{align*}  
with coproduct
\begin{align*}
\label{}
    \D{h_{i}}=&1\otimes h_{i} + h_{i} \otimes 1, & \D{e_{i}}=&e_{i}\otimes q^{\gamma_{i}}+ 1\otimes e_{i},
\end{align*}   
for all $i,j\in I$ and suitable elements $\gamma_{i}\in \h[[h]]$.
\end{lem}

\begin{pf}
After replacing the ordinary tensor product with the super-tensor product, the proof is identical to the proof of Proposition 3.1 of \cite{EK6}.  There are no new signs introduced.  For the most part, this is true because the arguments of the proof are based on the purely even Cartan subalgebra $\h$. 
\end{pf}

%   LEMMA  \gamma_{i}=d_{i}h_{i}

\begin{lem}\label{L:gamma}
$ \gamma_{i}=d_{i}h_{i}$
\end{lem}
\begin{pf}

By definition we have the natural projection $\borelhat_{+}\rightarrow\borel_{+}$.  Then the functoriality of the quantization implies that there is an epimorphism of Hopf superalgebras $U_{h}(\borelhat_{+}) \rightarrow U_{h}(\borel_{+})$.  Therefore $U_{h}(\borel_{+})$ is generated by $h_{i}, e_{i}$ satisfying the relations of Lemma \ref{L:UQESbhat} (and possibly other relations).  So it suffices to show that $ \gamma_{i}=d_{i}h_{i}$ in $U_{h}(\borel_{+})$.

Next we show that $U_{h}(\borel_{+}) \cong U_{-h}(\borel_{+})\qdo$.  From the definition of $\gl$ the Lie superbialgebra $\borel_{+}$ is self dual, i.e. $\borel_{+}\cong \borel_{+}^{*}$.  Again from functoriality we have that $U_{h}(\borel_{+}) \cong U_{h}(\borel_{+}^{*})$.  From Proposition \ref{P:glManinT} we have $\borel_{+}^{*}\cong\borel_{-}^{op}$.  Then equation (\ref{E:Ug+QDual}) and Theorem (\ref{T:twoQareIso})  imply that $U_{h}(\borel_{+})\qdo \cong U_{h}(\borel_{+}^{*op})$.  Substituting $\borel_{+}^{op}$ for $\borel_{+}$ we have $U_{h}(\borel_{+}^{op})\qdo \cong U_{h}(\borel_{+}^{*})$.  Finally from relation (\ref{R:ConditionDelta}) it follows that $U_{h}(\borel_{+}^{op}) \cong U_{-h}(\borel_{+})$ which implies that $U_{-h}(\borel_{+})\qdo \cong U_{h}(\borel_{+}^{*})$.  Thus, we have shown that $U_{h}(\borel_{+}) \cong U_{-h}(\borel_{+})\qdo$.  %Notice that in degree zero this isomorphism comes from the map $\h \rightarrow \h^{*}$ which corresponds to the form $2(,)$ of subsection \ref{SS:glbialg}.    

This isomorphism gives rise to the bilinear form $B: U_{h}(\borelp) \otimes U_{-h}(\borelp) \rightarrow \C((h))$ which satisfies the following conditions 
\begin{align}
\label{E:RelationsOfB1}
  B(xy,z)=& B(x \otimes y,\D z), & B(x, yz)=& B(\D x, y \otimes z) 
\end{align}
$$ B(q^{a},q^{b})=q^{-(a,b)}, a,b \in \h.$$
Let $a\in\h$ and $i\in I$.  Set $B_{i}=B(e_{i},e_{i})$, which is nonzero.   Using (\ref{E:RelationsOfB1}) we have 
\begin{align*}
\label{}
    B(e_{i},q^{a}e_{i})&=B(e_{i}\otimes q^{ \gamma_{i}} +1 \otimes e_{i}, q^{a}\otimes e_{i})   \\
    &  = B(e_{i}, q^{a})B(q^{ \gamma_{i}},e_{i}) + B(1,q^{a})B(e_{i},e_{i}) \\
    & = B_{i}
\end{align*}
since $B(e_{i},q^{a})=0$.  Similarly, we have $B(e_{i},q^{a}e_{i}q^{-a})= B(e_{i},q^{a}e_{i})B(q^{\gamma_{i}},q^{-a})$ implying
\begin{align}
\label{E:GammaB}
    B_{i}q^{(a,\gamma_{i})}=&B(e_{i},q^{a}e_{i}q^{-a}).  
\end{align}  
To complete the proof we need the following relation:
\begin{align}
\label{E:quantea}
   q^{a}e_{i}q^{-a}=&q^{\alpha_{i}(a)}e_{i} 
\end{align}
This relation is equivalent to $q^{h_{j}}e_{i}q^{-h_{j}}=q^{\alpha_{i}(h_{j})}e_{i}$ which follows from expanding $q=e^{h}$ and using the relation $[a,e_{i}]=\alpha_{i}(a)e_{i}$.
From (\ref{E:GammaB}) and (\ref{E:quantea}) we have 
$$B_{i}q^{(a,\gamma_{i})}=B(e_{i},q^{a}e_{i}q^{-a})=B(e_{i},q^{\alpha_{i}(a)}e_{i})=B_{i}q^{\alpha_{i}(a)}.$$
Thus, $(a, \gamma_{i})=\alpha_{i}(a)$, but $ \alpha_{i}(a)= d_{i}(a,h_{i})$, and so $\gamma_{i}=d_{i}h_{i}$, which completes the proof. 
\end{pf}

%    Subsection     The quantized enveloping superalgebra $\Uhb$

\subsection{The quantized universal enveloping superalgebra $\Uhb$}

In this subsection we show that there exist a bilinear form on $\Uhbhat$ such that $\Uhbhat$ modulo the kernel of the form is isomorphic to $\Uhb$.  

%  Thm   Form B

\begin{thm}\label{T:FormB}
There exists a unique bilinear form on $\Uhbhat$ which takes values in $\C((h))$ with the following properties
\begin{align*}
\label{}
  B(xy,z)=& B(x \otimes y,\D z), & B(x, yz)=& B(\D x, y \otimes z) 
\end{align*}
$$ B(q^{a},q^{b})=q^{-(a,b)}, a,b \in \h.$$
$$ B(e_{i},e_{j})=
	\begin{cases}
		(q_{i}-q_{i}^{-1})^{-1} & \text{if $i=j\neq m$,}\\
		1 & \text{if $i=j=m$.}\\
		0 & \text{otherwise}.
	\end{cases}$$
Moreover $\Uhb \cong \Up:=\Uhatp/Ker(B)$ as QUE superalgebras.	
\end{thm}
% PRoof of thm form B
\begin{pf}
The existence and uniqueness follows from the fact that the superalgebra generated by the $e_{i}$ is free. 

We will show that there is a nondegenerate bilinear form on $\Uhb$ with the same properties as $B$.  From the proof of Lemma \ref{L:gamma} we have that $U_{h}(\borel_{+}) \cong U_{-h}(\borel_{+})\qdo$.  But the even homomorphism $U_{-h}(\borel_{+})^{op} \rightarrow \Uhb$ given by conjugation by $q^{-\sum x_{i}^{2}/2}$, where $x_{i}$ is a orthonormal basis for $\h$, is a isomorphism.  Therefore we have a even isomorphism $\Uhb \cong \Uhb\qd$.  This isomorphism gives rise to the desired form on $\Uhb$.  

So the form $B$ is the pull back of the form on $\Uhb$.  Implying that the kernel of the form on $\Uhb$ is contained in the image of the kernel of $B$ under natural projection.

But the kernel of the form on $\Uhb$ is zero since the form is nondegenerate.   Thus we have $\Uhb \cong \Uhatp/Ker(B)$. 
\end{pf}

%  The kernel of $B$

\subsection{The kernel of $B$}  

In this subsection we show that $Ker(B)$ is generated by the quantum Serre-type relations (\ref{E:QserreA})-(\ref{E:QserreD}).   We first show that the quantum Serre-type relations are contained in $Ker(B)$.  To this end, we extend results of Lusztig \cite{L}.  The outline of this subsection is as follows.   
We start with the initial data: a free associative superalgebra $\Pf$ with unit and a Cartan matrix.  Using the Cartan matrix we define a twisted multiplication on $\Pf \otimes \Pf$ (see (\ref{E:tensorf})).  Then we prove that there is a unique form $\lform$ on $\Pf$ whose kernel contains the quantum Serre-type relations.  We end the subsection by showing that this implies that these relations are in $Ker(B)$.  
Intuitively, this construction is imposing the information of the Cartan matrix onto the twisted multiplication which in turn is imposing the relations on the kernel of $\lform$.

%Def of 'f
Let $q$ be an indeterminate.  Recall the definitions of Cartan data ($\Phi$, $(A,\tau)$,...) given at the beginning of this section.  Let $\Pf$ be the free associative $\C(q)$-superalgebra with 1 generated by $\theta_{i}$, for $i \in I$, where the parity is $\p 0$ for all generators except for $\theta_{i}, \: i\in \tau$ which has parity $\p1$. 

%Def of f_{v}
 
 For any $\nu=\sum_{i}\nu_{i}i \in \N[I]$,  let $\Pf_{\nu}$ be the $\C(q)$-subspace of $\Pf$ spanned by the monomials $\theta_{i_{1}}\theta_{i_{2}}\dots\theta_{i_{k}}$ so that for each $i \in I$, the number of times $i$ appears in the sequence $i_{1}, i_{2}, \dots, i_{k}$ is equal to $\nu_{i}$.  Notice that $\Pf=\oplus_{\nu}\Pf_{\nu}$.  We say $x\in \Pf$ is homogeneous if $x\in \Pf_{\nu}$, for such an $x$, set $|x|=\nu$.  For homogeneous $x, x' \in \Pf$, let 
 \begin{equation}
\label{E:Def<,>}
<|x|,|x'|>:=<\sum_{i} d_{i}\nu_{i}h_{i},\sum_{j} \nu'_{j}\alpha_{j}>=\sum_{i,j}d_{i}\nu_{i}\nu'_{j}\alpha_{j}(h_{i}),
\end{equation}
  where $|x|=\sum \nu_{i}i$ and $|x'|=\sum \nu'_{j}j$.  Note that $<|x|,|x'|>=(\sum \nu_{i}h_{i},\sum \nu'_{j}h_{j})$ where $(,)$ is the super-symmetric bilinear form on $\g$ (see \ref{SS:glbialg}). % To be consistent with \cite{L} we use  (\ref{E:Def<,>}) to define $<,>$. %$<|\theta_i|,|\theta_j|>=<d_{i}h_{i},\alpha_{j}>=d_{i}a_{ij}=d_{j}a_{ji}=<|\theta_{j}|,|\theta_{i}|>$. 
 
 We make $\Pf \otimes \Pf$ into an a superalgebra with the following multiplication:
\begin{equation}
\label{E:tensorf}
(x_{1}\otimes x_{2})(y_{1 } \otimes y_{2}) = (-1)^{\p x_{2} \p y_{1} }q^{<|x_{2}|, |y_{1}|>}x_{1}y_{1}\otimes x_{2}y_{2}.
\end{equation}
where $x_{1},x_{2},y_{1},y_{2}\in \Pf $ are homogeneous.
%In particular, $$(1 \otimes \theta_{i})(\theta_{j}\otimes 1)=(-1)^{\p  \theta_{i} \p \theta_{j} } q^{<d_{i}h_{i},\alpha_{j}>}(\theta_{j}\otimes \theta_{i})=(-1)^{\p  \theta_{i} \p \theta_{j} } q^{d_{i}a_{ij}}(\theta_{j} \otimes \theta_{i}).$$

%Def of r

Let $r:\Pf \rightarrow \Pf \otimes \Pf $ be the superalgebra map defined by $r(\theta_{i})=\theta_{i}\otimes 1 + 1\otimes \theta_{i}$.

\begin{prop}
\label{P:formf}
There is a unique bilinear form $\lform$ on $\Pf$ with values in $\C(q)$ such that $\lform(1,1)=1$ and 
\begin{enumerate}

  \item $
	 \lform(\theta_{i},\theta_{j})=
	\begin{cases}
		(q_{i}-q_{i}^{-1})^{-1} & \text{if $i=j\neq m$,}\\
		1 & \text{if $i=j=m$.}\\
		0 & \text{otherwise},
	\end{cases}
$

\item $ \lform(x,yz )=\lform(r(x),y \otimes z)$ for all $x,y,z \in \Pf$, \label{Pitem:r1}
\item $ \lform(xy,z )=\lform(x \otimes y, r(z))$ for all $x,y,z \in \Pf$, \label{Pitem:r2}

\end{enumerate}
where the bilinear form on $\Pf \otimes \Pf$ (also denoted by $\lform$) is given by 
\begin{equation}
\label{E:bilinearform}
\lform(x_{1}\otimes x_{2},y_{1 } \otimes y_{2}) = (-1)^{\p x_{2} \p y_{1} }\lform(x_{1},y_{1})\lform(x_{2},y_{2}).
\end{equation} 
\end{prop}

%proof of (,)

\begin{pf}

The proof is similar to the proof of Proposition 1.2.3 in \cite{L}.   Here we define $\lform$ and refer the reader to \cite{L} for the rest of the proof.  First, we define a superalgebra structure on $\Pf^{*}$. 

 For any $\nu, \nu^{\prime} \in \N[I]$, composing the map $r |_{\Pf_{\nu+\nu^{\prime}}}: \Pf_{\nu+\nu^{\prime}}\rightarrow \Pf \otimes \Pf$ with the projection $\Pf\otimes \Pf \rightarrow \Pf_{\nu}\otimes \Pf_{\nu^{\prime}}$, we have the linear map 
$\Pf_{\nu+\nu^{\prime}}\rightarrow \Pf_{\nu}\otimes \Pf_{\nu^{\prime}}$.  Taking the dual, we obtain linear maps $\Pf_{\nu}^{*}\otimes \Pf_{\nu^{\prime}}^{*}\rightarrow \Pf_{\nu+\nu^{\prime}}^{*}$.  This defines an associative superalgebra structure on $\Pf^{*}$.
For each $i\in I$, let $\xi_{i}\Pf_{i}^{*}$ be given by 
$$
	 \xi_{i}(\theta_{i})=
	\begin{cases}
		(q_{i}-q_{i}^{-1})^{-1} & \text{if $i\neq m$,}\\
		1 & \text{if $i=m$.}
	\end{cases}
$$.

Let $\phi:\Pf \rightarrow \oplus_{\nu}\Pf^{*}_{\nu}$ be the unique superalgebra homomorphism preserving 1, such that $\phi(\theta_{i})=\xi_{i}$ for all $i$.  For homogeneous $x,y \in \Pf$, set $\lform(x,y)=(-1)^{\p x \p y}\phi(y)(x)$.  Now (\ref{Pitem:r1}) follows as $\phi$ is an algebra homomorphism.  From the definition of $\phi$ we have 
\begin{equation}
\label{E:MustBeHomo}
\lform(x,y)=0 \text{ unless } |x|=|y|
\end{equation}
for homogeneous $x,y \in \Pf$.
  
After putting in appropriate signs coming from (\ref{E:tensorf}) and (\ref{E:bilinearform}), the proof of (\ref{Pitem:r2}) follows as in \cite{L}.  
\end{pf}

%def  I and f

Let $Ker(\lform)$ be the kernel of the form $\lform$, then $Ker(\lform)$ is a homogeneous ideal of $\Pf$.  Let $\f=\Pf / Ker(\lform)$.  From (\ref{E:MustBeHomo}),  the decomposition 
$\Pf = \oplus_{\nu}\Pf_{\nu}$ gives a direct sum decomposition of $\f = \oplus_{\nu}\f_{\nu}$, where $\f_{\nu}$ is the image of $\Pf_{\nu}$ under the projection $\Pf\rightarrow \f$.

%prop serre rel in ker of (,)

\begin{prop}\label{P:QSinf}

The relations (\ref{E:QserreA})-(\ref{E:QserreD}) with $e$ replaced by $\theta$ hold in the superalgebra $\f$.  In particular, 
\begin{equation}
  \label{R:even0}
   \theta_{m}^{2}= 0 \quad \text{ if $\tau=\{m\}$,}
\end{equation}
\begin{multline}
\label{E:QserreCtheta}
\theta_{m}\theta_{m-1}\theta_{m}\theta_{m+1}+\theta_{m}\theta_{m+1}\theta_{m}\theta_{m-1}+\theta_{m-1}\theta_{m}\theta_{m+1}\theta_{m}
+\theta_{m+1}\theta_{m}\theta_{m-1}\theta_{m}\\-(q+q^{-1})\theta_{m}\theta_{m-1}\theta_{m+1}\theta_{m}=0  \quad \text{ if $ m-1,m, m+1\in I $ and $  a_{mm}=0$,}
\end{multline}
\begin{multline}
\label{E:QserreDtheta}
\theta_{m-1}\theta_{m}^{3}- (q+q^{-1}-1)\theta_{m}\theta_{m-1}\theta_{m}^{2}
-(q+q^{-1}-1)\theta_{m}^{2}\theta_{m-1}\theta_{m} + \theta_{m}^{3}\theta_{m-1}=0\\ \quad \text{ if the Cartan Matrix $A$ is of type B, $\tau = \{m\}$ and $s=m$}.
\end{multline}
\end{prop}

%proof that serre rel are in ker
\begin{pf}
The relations (\ref{E:QserreA2}) and (\ref{E:QserreB}) $e$ replaced with $\theta$ are the normal quantum Serre relations and follow from Proposition 1.4.3 in \cite{L}.   From  (\ref{E:MustBeHomo}), we have relation (\ref{R:even0}) holds if $\lform(\theta^{2}_{m},\theta^{2}_{m})=0$, which follow immediately from Proposition \ref{P:formf}, part (\ref{Pitem:r1}).  

We will now show the relation (\ref{E:QserreCtheta}) holds when $ m-1,m, m+1\in I $ and $  a_{mm}=0$.   In this case we have 
\begin{equation}
\label{E:FormulaForA}
a_{ij}=(1+(-1)^{\delta_{i,m}})\delta_{i,j}-(-1)^{\delta_{i,m}}\delta_{i,j-1}-\delta_{i,j+1}
\end{equation}
for $i,j\in \{m-1,m,m+1\}$.  We also have $d_{m-1}=d_{m}=1$ and $d_{m+1}=-1$.

Let $l$ be the left side of relation (\ref{E:QserreCtheta}).
To show that relation (\ref{E:QserreCtheta}) holds it is enough to show $\lform(x,l)=0$ for all $x\in \Pf_{1(m-1)+2(m)+1(m-1)}$.  By relation (\ref{R:even0}) the vector space $\f_{1(m-1)+2(m)+1(m-1)}$ is generated by 
$$\theta_{m}\theta_{m-1}\theta_{m}\theta_{m+1}, \theta_{m}\theta_{m+1}\theta_{m}\theta_{m-1},$$ $$ \theta_{m-1}\theta_{m}\theta_{m+1}\theta_{m}, \theta_{m+1}\theta_{m}\theta_{m-1}\theta_{m},  \theta_{m}\theta_{m-1}\theta_{m+1}\theta_{m}.$$
%$\theta_{m}\theta_{m-1}\theta_{m}\theta_{m+1},$ $ \theta_{m}\theta_{m+1}\theta_{m}\theta_{m-1},$ $\theta_{m-1}\theta_{m}\theta_{m+1}\theta_{m}$,  $\theta_{m+1}\theta_{m}\theta_{m-1}\theta_{m}$,  $ \theta_{m}\theta_{m-1}\theta_{m+1}\theta_{m}$.   
Therefore it suffices to check that $\lform(x,l)=0$, when $x$ is any of the above generators.  We will check this condition for $\theta_{m +1 }\theta_{m  }\theta_{m -1 }\theta_{m  }$, the others follow similarly.

Let $c_{i}=(q_{i}-q_{i}^{-1})^{-1}$.  From (\ref{E:tensorf}) we have
\begin{equation}
\label{E:tensorOfThetas}
(1 \otimes \theta_{i})(\theta_{j}\otimes 1)=(-1)^{\p  \theta_{i} \p \theta_{j} } q^{d_{i}a_{ij}}(\theta_{j} \otimes \theta_{i}).
\end{equation}

We use (\ref{E:bilinearform}), (\ref{E:FormulaForA}), (\ref{E:tensorOfThetas}) and Proposition \ref{P:formf}, part (\ref{Pitem:r2}) to make the following calculations:
\begin{align*}
\label{}
    a_{1}:=&\lform(\theta_{m +1 }\theta_{m  }\theta_{m -1 }\theta_{m  },\theta_{m  }\theta_{m -1 }\theta_{m  }\theta_{m +1 })\\
    =&\lform(\theta_{m +1 }\theta_{m  } \otimes \theta_{m -1 }\theta_{m  },r(\theta_{m  }\theta_{m -1 })r(\theta_{m  }\theta_{m +1 }))   \\
    =&-q\lform(\theta_{m +1 }\theta_{m  },\theta_{m  }\theta_{m +1 })\lform(\theta_{m -1 }\theta_{m  },\theta_{m -1 }\theta_{m  })\\
     & \hspace{30pt}-(-qq^{-1})\lform(\theta_{m +1 }\theta_{m  },\theta_{m  }\theta_{m +1 })\lform(\theta_{m -1 }\theta_{m  },\theta_{m  }\theta_{m-1  })  \\
    =&-q(qc_{m+1})(c_{m-1})-(-qq^{-1})(qc_{m+1})(q^{-1}c_{m-1})\\
    =&-q^{2}c_{m+1}c_{m-1}+c_{m+1}c_{m-1}.  
\end{align*}
Similarly we have
\begin{align*}
\label{}
    a_{2}:=&\lform(\theta_{m +1 }\theta_{m  }\theta_{m -1 }\theta_{m  },\theta_{m  }\theta_{m+1  }\theta_{m }\theta_{m-1  })=0   \\
    a_{3}:=&\lform(\theta_{m +1 }\theta_{m  }\theta_{m -1 }\theta_{m  },\theta_{m -1 }\theta_{m  }\theta_{m +1 }\theta_{m  })=0   \\
    a_{4}:=&\lform(\theta_{m +1 }\theta_{m  }\theta_{m -1 }\theta_{m  },\theta_{m +1 }\theta_{m  }\theta_{m -1 }\theta_{m  })=-c_{m+1}c_{m-1}+q^{-2}c_{m+1}c_{m-1}   \\
    a_{5}:=&\lform(\theta_{m +1 }\theta_{m  }\theta_{m -1 }\theta_{m  },\theta_{m  }\theta_{m-1  }\theta_{m +1 }\theta_{m  })=-qc_{m+1}c_{m-1}+q^{-1}c_{m+1}c_{m-1}  
\end{align*}
So\begin{multline*}
\label{}
    \lform(\theta_{m +1 }\theta_{m  }\theta_{m -1 }\theta_{m  },l)=a_{1}+a_{2}+a_{3}+a_{4}-(q+q^{-1})a_{5} \\
    =[-q^{2}+1-1+q^{-2}-(q+q^{-1})(-q+q^{-1})]c_{m+1}c_{m-1}=0.  
\end{multline*}
It is not hard to follow the above computation and show that (\ref{E:QserreDtheta}) holds and so the Proposition follows.
\end{pf}

%Next we define the Drinfeld-Jimbo type superalgebra $\Ug$ over the indeterminate $q$ (for more details see \cite{FLV},\cite{KTol}).  Note $\Ug$ is an algebraic construction where $\DJg$ is a ``topological'' construction.  Let $\Ug$ be the $\C(q)$-superalgebra generated by $q^{h_{i}}$, $e_{i} \text{ and } f_{i} $ for $i \in I$  modulo relations (\ref{E:DJglRelation2})-(\ref{E:QserreD}) and 

One can continue to follow \cite{L} and show that the Drinfeld-Jimbo type $\C(q)$-superalgebra (see \cite{FLV,KTol}) can be recovered from $\f$.  This result is not essential for our purposes here.  However, in order to shed some light on the larger picture we will now state the results without proof. 

Let $\pu$ be the $\C(q)$-superalgebra generated by $q^{h_{i}}$, $e_{i} \text{ and } f_{i} $ for $i \in I$ modulo the relations (\ref{E:DJglRelation2}) and 
\begin{align}
\label{R:DJq1}
    q^{h_{i}}q^{h_{j}} =&  q^{h_{i}+h_{j}}, 
&   q^{h_{i}}e_{j} &=   q^{a_{ij}}e_{j}q^{h_{i}}, & q^{h_{i}}f_{j}& =  q^{-a_{ij}}f_{j}q^{h_{i}}.
\end{align}  Let $\cu$ be the associative $\C(q)$-superalgebra $\pu$ modulo the following relations:
for any relation $g(\theta_{i})\in  Ker(\lform)$ we have $g(e_{i})=0$ and $g(f_{i})=0$ in $\cu$.  
Let $\cu^{0} $ ($\pu^{+}_{0}$) be the sub-superalgebra of $\pu$ generated by $q^{h_{i}}\: (i\in I)$ (resp. $q^{h_{i}}, e_{i}\: (i\in I)$).  
Let $\f \rightarrow \cu \; (x\mapsto x^{+})$ 
 and $\f \rightarrow \cu \; (x\mapsto x^{-})$ 
  be the homomorphism such that $e_{i}=\theta_{i}^{+}$ and $f_{i}=\theta_{i}^{-}$ for all $i \in I$.  As in \cite{L}, one can show that 
 \begin{align*}
\label{}
    \f \otimes \cu^{0} \otimes \f \rightarrow  \cu  &\text{ given by } u \otimes q^{a} \otimes w \mapsto u^{-}q^{}w^{+} \\
    \cu^{0}\otimes \Pf \rightarrow  \pu^{+}_{0} &\text{ given by } q^{a}\otimes x \mapsto q^{a}x^{+}
\end{align*}
 are isomorphisms of vector spaces.  From the first isomorphism and Proposition \ref{P:QSinf} it follows that the superalgebra $\cu$ is isomorphic to the D-J type $\C(q)$-superalgebra (which is the superalgebra $\pu$ modulo (\ref{E:QserreA})-(\ref{E:QserreD})).  

Now we are ready to prove the following theorem.

\begin{thm}\label{T:KerBContainQSerre}
The quantum Serre-type relations (\ref{E:QserreA})-(\ref{E:QserreD}) are contained in $Ker(B)$.
\end{thm}
\begin{pf}
Recall the superalgebra $\Uhatp$ of Theorem \ref{T:GRofUhbhat}.  By setting $q$ to $e^{h/2}$ one can obtain an injective superalgebra morphism $\pu^{+}_{0} \rightarrow \Uhatp$.  Then the composition $$\Pf \hookrightarrow  \cu^{0}\otimes \Pf \rightarrow \pu^{+}_{0} \rightarrow \Uhatp$$ 
is injective.  The form $\lform:\Pf\otimes \Pf\rightarrow \Pf$ of Proposition \ref{P:formf} corresponds (under the above composition) to the form $B$.  Therefore, Proposition \ref{P:QSinf} implies that the quantum Serre-type relations are contained in $Ker(B)$
\end{pf}

\begin{cor}
\label{T:KerBQSerre}
$Ker(B)$ is generated by the quantum Serre-type relations (\ref{E:QserreA})-(\ref{E:QserreD}).
\end{cor}

\begin{pf}
Let $\Up=\Uhatp/Ker(B)$.  By construction the superalgebra $\Uhb$ is isomorphic as a vector space to $U(\borelp)[[h]]$, implying $\Up\cong U(\borel_{+})[[h]]$.  Combining this observation with Theorem \ref{T:KerBContainQSerre} and the fact that $\borelp$ is the quotient of $\borelhat_{+}$ by the classical super Serre-type relations (\ref{E:ClassicSerre}) we have that the $Ker(B)$ is generated by the quantum Serre-type relations.  
\end{pf}

%    Subsection     The quantized enveloping superalgebra U_{h}(\gl)

\subsection{Generators and relations for $U_{h}(\g)$}
\begin{thm}\label{T:UhIsoDJ}
Let $\g$ be a Lie superalgebra of type A-G. The QUE superalgebra $U_{h}(\g)$ is isomorphic to the quotient of the double $D(\Up)$ by the ideal generated by the identification of $\h \subset \Up$ and $\h^{*} \subset \Upqd$, i.e. the Etingof-Kazhdan quantization $U_{h}(\g)$ is isomorphic to the Drinfeld-Jimbo type superalgebra $\DJg$ (see \S \ref{SS:GenRel}). 
\end{thm}
\begin{pf}
Recall from \S \ref{SS:glbialg} that the Lie superbialgebra structure of $\g$ comes from identifying $\h$ and $\h^{*}$ in $\g \oplus \h= \borel_{+}\oplus \borel_{+}^{*}$.  Also since the quantization commutes with the double we have 
$$\OUh(D(\borel_{+})) \cong D(\Uhb)=\Uhb \otimes \Uhb\qdo.$$ 
Therefore, we have $U_{h}(\g)$ is isomorphic to $D(\Uhb)=\Uhb \otimes \Uhb\qdo$ modulo the the ideal generated by the identification of $\h \subset \Uhb$ and $\h^{*} \subset \Uhb\qdo$.  But from Theorem \ref{T:FormB} we have that $D(\Uhb)\cong D(\Up)$ and then Corollary \ref{T:KerBQSerre} implies result.
\end{pf}

\section{A theorem of Drinfeld's}\label{S:ModuleCat}
Recall the definition of $\kzt$ and $\DJg$ given in \S\ref{SS:BQBAU} and \S\ref{SS:GenRel} respectively.  Here we use all the results of this paper to show that the categories of topologically free modules over $\kzt$ and $\DJg$ are braided tensor equivalent.  We do this in two steps: (1) we show that $\Uh$ and $\kzt$ have equivalent module categories, (2) we use the fact the that $\DJg$ and $\Uh$ are isomorphic to prove the desired result.   For more on braided tensor categories see \cite{ES,Kas}. 

\subsection{The E-K quantization $\Uh$ and $\kzt$}
In this subsection we show that $\Uh$ is the twist of $\kzt$ by $J$.  To this end we recall the following definitions.  

Let $(A,\Delta,\epsilon, \Phi,R)$ be a quasitriangular quasi-superbialgebra (see \S\ref{SS:BQBAU}.)  An invertible element $J\in A\otimes A$ is a \emph{gauge transformation on $A$} if 
$$(\epsilon \otimes id)(J)=(id \otimes \epsilon)=1.$$
Using a gauge transformation $J$ on $A$, one can construct a new quasitriangular quasi-superbialgebra $A_{J}$ with coproduct $\Delta_{J}$, R-matrix $R$ and associator $\Phi_{J}$ defined by
$$\Delta_{J}= J^{-1}\Delta J, \:\: R_{J}=(J^{op})^{-1}RJ,$$
$$  \Phi_{J}= J^{-1}_{23}(id \otimes \Delta)(J^{-1})\Phi(\Delta \otimes id)(J)J_{12}.$$

As is the case of quasitriangular (quasi-)bialgebra, the category of modules over a quasitriangular (quasi-)superbialgebra is a braided tensor category.

\begin{thm}\label{T:EquivOfqqbialg}
Let $A$ and $A'$ be a quasitriangular quasi-superbialgebra.  Suppose that $J$ is a gauge transformation on $A'$ and $\alpha:A\rightarrow A_{J}'$ is an isomorphism of quasitriangular quasi-superbialgebra then $\alpha$ induces a equivalence between the braided tensor categories $A'\text{-}Mod$ and $A\text{-}Mod.$
\end{thm}
\begin{pf}
Let $\alpha^{*}:A'\text{-}Mod \rightarrow A\text{-}Mod$ be the functor defined as follows.  On objects, the functor $\alpha^{*}$ is defined by sending the module $W$ to the same underlying vector space with the action given via the isomorphism $\alpha$.  For any morphism $f: W\rightarrow X$ in $A'\text{-}Mod$ let $\alpha^{*}(f)$ be the image of $f$ under the isomorphism
$$\Hom_{A'}(W,X)\cong \Hom_{A}(W,X).$$
A standard categorical argument shows that this functor is an equivalence of braided tensor categories (see \S XV.3 of \cite{Kas}).
\end{pf}

Let $\g$ be a Lie superalgebra of type A-G.  Recall from \S \ref{SS:glbialg} that $\g$ has a unique non-degenerate supersymmetric invariant bilinear form.  Let $t$ be the corresponding even invariant super-symmetric element of $\g \otimes \g$.  Let $\OJ$ be the element of $\U[[h]]^{\otimes 2}$ defined in (\ref{E:defJ}).  By definition of the coproduct and R-matrix of $\Uh$ (see \S\ref{SS:Qofthedouble}) we have that $\Uh=(\kzt)_{\OJ}$.

\subsection{Main theorem}
Let $X$ be a topological (quasi) Hopf superalgebra and let $ X \text{-}Mod_{fr}$ of topologically free $X$-modules of finite rank (see \S\ref{SS:TopoFreeMod}). 
%We now give an equivalence of the braided tensor categories $\kzt\text{-}Mod_{fr}$ and $\DJg\text{-}Mod_{fr}$. 
The following theorem was first due to Drinfeld \cite{D5} in the case of semi-simple Lie algebras.

\begin{thm}\label{T:Drinfeldsuper}
The braided tensor categories $\kzt\text{-}Mod_{fr}$ and $\DJg\text{-}Mod_{fr}$ are equivalent. 
\end{thm} 
\begin{pf}
As mentioned at the end of the last subsection $\Uh=(\kzt)_{\OJ}$.  Combining this fact with 
Theorem \ref{T:UhIsoDJ} we have that there exists an isomorphism of quasitriangular quasi-superbialgebra
$$\alpha: \DJg \rightarrow (\kzt)_{\OJ}.$$ 
Now as a consequence of Theorem \ref{T:EquivOfqqbialg} we have that the categories $\kzt\text{-}Mod_{fr}$ and $\DJg\text{-}Mod_{fr}$ are braided tensor equivalent.
\end{pf}

\begin{rem}
Drinfeld's proof of Theorem \ref{T:Drinfeldsuper} in the case of semi-simple Lie algebras uses deformation theoretic arguments to show the existence of $\alpha$.  Our proof constructs the isomorphism $\alpha$ explicitly.
\end{rem}

\end{document}